\documentclass[journal]{IEEEtran}

\usepackage{amsmath}
\usepackage{dsfont}
\usepackage{graphicx}
\usepackage{hyperref}
\usepackage{bm}
\usepackage{tcolorbox}

\newcommand{\matr}[1]{\mathbf{#1}}

\usepackage[ruled,norelsize]{algorithm2e}
\makeatletter
\newcommand{\removelatexerror}{\let\@latex@error\@gobble}
\makeatother
\begin{document}

\title{Nonlinear control of networked dynamical systems}

\author{Megan~Morrison 
        and~J.~Nathan~Kutz% <-this % stops a space
\thanks{M. Morrison and J.N. Kutz are with the Department
of Applied Mathematics, University of Washington, Seattle,
WA, 98195 USA e-mail: mmtree@uw.edu and kutz@uw.edu.}}

% The paper headers
\markboth{IEEE TRANSACTIONS ON NETWORK SCIENCE AND ENGINEERING}%
{Shell \MakeLowercase{\textit{et al.}}: Bare Demo of IEEEtran.cls for IEEE Journals}
%
% make the title area
\maketitle
%
% %%%%%%%%%%%% to do list %%%%%%%%%
% \listoffixmes
% \vspace{1cm}
% %%%%%%%%%%%%%%%%%%%%%%%%%%%%%%%%%
%
% As a general rule, do not put math, special symbols or citations
% in the abstract or keywords.
\begin{abstract}
We develop a principled mathematical framework for controlling nonlinear, networked dynamical systems. Our method integrates dimensionality reduction, bifurcation theory and emerging model discovery tools to find low-dimensional subspaces where feed-forward control can be used to manipulate a system to a desired outcome. The method leverages the fact that many high-dimensional networked systems have many fixed points, allowing for the computation of control signals that will move the system between any pair of fixed points. The {\em sparse identification of nonlinear dynamics} (SINDy) algorithm is used to fit a nonlinear dynamical system to the evolution on the dominant, low-rank subspace. This then allows us to use bifurcation theory to find collections of constant control signals that will produce the desired objective path for a prescribed outcome. Specifically, we can destabilize a given fixed point while making the target fixed point an attractor. The discovered control signals can be easily projected back to the original high-dimensional state and control space. We illustrate our nonlinear control procedure on established bistable, low-dimensional biological systems, showing how control signals are found that generate switches between the fixed points. We then demonstrate our control procedure for high-dimensional systems on random high-dimensional networks and Hopfield memory networks.
\end{abstract}
%
% Note that keywords are not normally used for peerreview papers.

\IEEEpeerreviewmaketitle

\section{Introduction}
\IEEEPARstart{N}{etworked} dynamical systems are ubiquitous across the engineering, physical, biological and social sciences. They are often characterized by a high-dimensional state space and nonlinearity, making them exceptionally difficult to characterize and control.  Indeed, it is typical that the connectivity is so complex that the functionality, control and robustness of the network of interest is impossible to characterize using standard mathematical methods.  Moreover, with few exceptions, underlying nonlinearities impair our ability to construct analytically tractable solutions, forcing one to rely on experiments and/or modern high-performance computing to study a given system. Unlike engineered systems that are constructed to be both measurable and controllable, such emergent systems can be difficult to measure and have restricted avenues of control. However, advances over the past decade have revealed a critical observation, that meaningful input/output of signals in high-dimensional networks are often encoded in low-dimensional patterns of dynamic activity~\cite{jones_natural_2007,rabinovich_transient_2008,brunton_data-driven_2019,kato_global_2015, marvel_continuous-time_2011,fieseler_unsupervised_2020,delahunt_biological_2018}.  We show that such low-dimensional patterns of activity can be exploited in order to develop principled techniques for a feed-forward architecture for establishing control of high-dimensional, nonlinear networked dynamical systems.  

The potential applications of a control framework for networked dynamical systems are extensive.  Neuroscience is an especially relevant example where networks of neurons interact to encode and process input stimulus and behavioral responses.  Recent observations in a variety of organisms, from the nematode {\em C. elegans}~\cite{kunert_low-dimensional_2014,kato_global_2015,fieseler_unsupervised_2020} to insect olfactory processing~\cite{jones_natural_2007,rabinovich_transient_2008,shlizerman_data-driven_2014,delahunt_biological_2018}, shows that the underlying encodings and control are fundamentally low-dimensional.  Network models are also common in attempts to understand the formation and retrieval of memories, such as proposed in the Hopfield model where each memory is a fixed point in the high-dimensional, networked dynamical system~\cite{hopfield_neural_1982, morrison_preventing_2017}.  Indeed, it is known that the nervous systems carries out an impressive feat of dimensionality reduction when it encodes behavior, collapsing the high-dimensional representation of the stimulus environment into the much lower representations for decision making and motor command.  Practical emerging technologies, such as {\em deep brain stimulation} (DBS), aim to leverage such control protocols to restore patients to their original functional capabilities.   Applications extend well beyond neuroscience, with the potential of the method to impact disease modeling~\cite{kermack_contribution_1927, hethcote_mathematics_2000}, social~\cite{morrison_community_2019} and financial networks, powergrid networks~\cite{dylewsky_engineering_2019}, and ecosystems, for instance.  

Characterization is only the first step in understanding networked dynamical systems.  A principled quantification of the low-dimensional patterns of dynamic activity can help lead to control protocols for manipulating the system into a desired outcome.  An extensive body of literature exists on the analysis and control of nonlinear systems \cite{isidori_nonlinear_1990}.  Unlike many linear control models, where controllability and observability can be explicitly computed and guaranteed, nonlinear control remains challenging, especially in networked settings.  Nonlinear networked dynamical systems can contain many fixed points, limit cycles and strange attractors, all of which make the development of principled control models difficult.  The manifestation of these various phenomenon must be addressed in any practical control paradigm.  On the other hand, one can use the existence of such rich dynamical structures to allow the network itself, under suitable manipulation, to evolve to a desired state of behavior.  Thus nonlinearity can exploit a much broader class of dynamics and function than linear models.
\par
We integrate methods of dimensionality reduction and data-driven discovery of dynamics to construct principled methods for controlling nonlinear, networked dynamical systems.  Specifically, we develop feed-forward control techniques to interpret and regulate the dynamics of such systems by leveraging dominant, low-dimensional subspaces on which the dynamics evolves. Our mathematical architecture generates a set of actuation signals that, when applied, are able to control the original high-dimensional system. Using bifurcation theory, we find  collections of feed-forward control signals that will force convergence to desired objective states, allowing us to move the system from one fixed point of the system to another in a principled manner. Specifially, we can destabilize a given fixed point by making it undergo a saddle node or Hopf bifurcation, while simultaneously making the target fixed point an attractor.  This creates a pathway with the feed-forward signals from one fixed point to another.  We first demonstrate our nonlinear control procedure on established bistable, low-dimensional biological systems showing how control signals are found that generate switches between attractors. We then show how random high-dimensional networks and Hopfield memory networks can be reliably controlled by discovering low-dimensional subspaces which characterize their dynamic evolution.  Our algorithmic procedure is the first of its kind to provide a principled mathematical architecture that simultaneously leverages mode discovery, dimensionality reduction, and bifurcation theory for controlling networked dynamical systems.

The paper is outlined as follows:  Sec.~II provides a brief overview of the background material necessary for constructing our control framework.  Section~III provides an analysis of feed-forward control applied to low-dimensional nonlinear dynamical systems.  This highlights the basic mathematical architecture that is used in Sec.~IV and V for network control applications in both low- and high-dimensional systems respectively.  A number of practical applications are considered in Sec.~VI, including the Hopfield memory model where we show how control can be used to transition between fixed points, or memories, in the network.  The paper is concluded in Sec.~VII with a discussion of our results.

\section{Background}

\subsection{Dimensionality Reduction}
It is typically observed that high-dimensional dynamical systems manifest behavior on low-dimensional manifolds~\cite{jones_natural_2007,rabinovich_transient_2008,brunton_data-driven_2019,kato_global_2015, marvel_continuous-time_2011,fieseler_unsupervised_2020,delahunt_biological_2018}.  Indeed, low-dimensional structures can be exploited for characterizing pattern forming systems~\cite{cross_pattern_1993} and reduced order models~\cite{benner_survey_2015}. {\em Principal component analysis} (PCA) is a linear dimensionality reduction technique based upon the {\em singular value decomposition} (SVD) that can extract dominant correlated features in complex, high-dimensional data \cite{jolliffe_principal_2002, kutz_data-driven_2013, brunton_data-driven_2019}, thus producing a coordinate system (subspace) on which our networked dynamics of interest can be projected.  Let ${\bf X} \in \mathds{R}^{n \times m}$ represent timeseries data collected from an $n$ dimensional system for $m$ timepoints. The SVD of this matrix produces the matrix decomposition~\cite{trefethen_numerical_1997, kutz_data-driven_2013, brunton_data-driven_2019}
\begin{align}
    {\bf X} = {\bf U}{\bf S}{\bf V}^*
\end{align}
were ${\bf U}$ denotes the dominant correlated spatial structures of the $n$-dimensional system, ${\bf S}$ is a diagonal matrix whose singular values characterize an ordered ranking of the correlations, and  ${\bf V}$ represent the projection of the modes into the temporal dimension.  Both ${\bf U}$ and ${\bf V}$ are unitary matrices with orthonormal columns.
\begin{figure}[t]
    \centering
    \includegraphics[width = 1.0\linewidth]{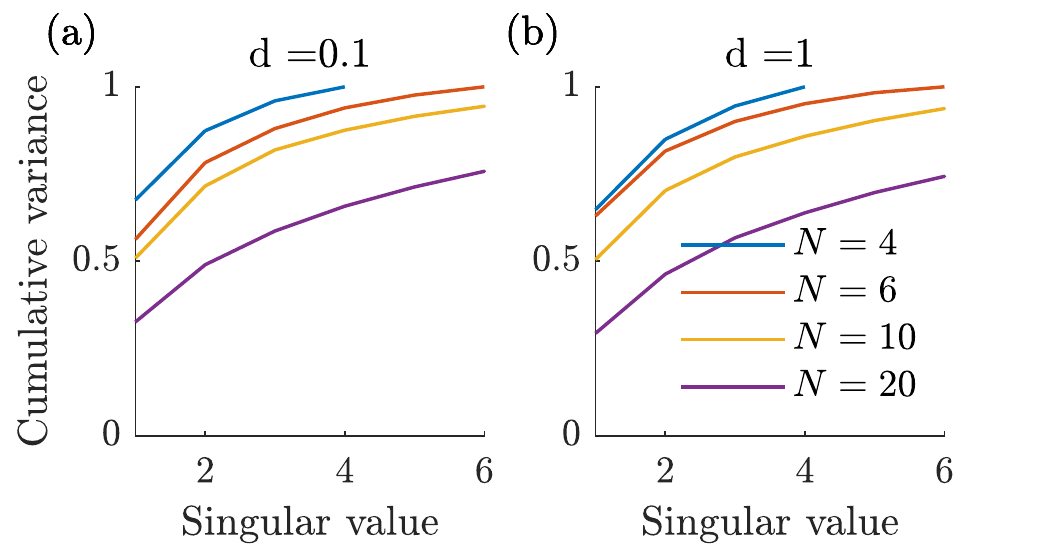}
    \caption{Average cumulative variance captured by initial singular values in random dynamical systems of increasing size $N = 4, 6, 10, 20$. (a) Dynamical systems consisting of a
    low percentage of possible term combinations, $d = 0.1$ (b) All possible term combinations included in the dynamical system, $d = 1.0$. Variance averaged over 100 trials}
    \label{fig:dimension}
\end{figure}
A low-dimensional, $r$-rank system can be optimally approximated in an $\ell_2$-sense using the first $r$ columns of each matrix: ${\bf X} = \hat{\bf U} \hat{\bf S} \hat{\bf V}^*$ where $\hat{{\bf U}} \in \mathds{R}^{n\times r}$, $\hat{{\bf S}} \in \mathds{R}^{r\times r}$, and $\hat{{\bf V}}^* \in \mathds{R}^{r\times m}$.  Such low-rank subspaces are exploited for building reduced order models that approximate the high-dimensional system~\cite{benner_survey_2015,brunton_data-driven_2019}.  It is also exploited in what follows since random networks generically manifest low-dimensional behavior, as shown in Fig.~\ref{fig:dimension}  where the cumulative variance contained in the initial modes of a randomly generated dynamical system is plotted.

\subsection{Sparse Identification of Nonlinear Dynamics - SINDy}
SINDy is a data-driven approach to find the sparse dynamics driving a dynamical system purely from timeseries data \cite{brunton_discovering_2016}.
If $\dot{x} = f(x)$ is the unknown true dynamics of a system from which we capture timeseries measurements ${\bf X}$, the SINDy algorithm seeks to approximate $f(x)$ by a generalized linear model in a set of candidate basis functions $\theta_k({x})$
	\begin{equation}
	{f}({x}) \approx \sum_{k=1}^p \theta_k({x})\xi_k = \boldsymbol{\Theta}({x})\mathbf{\Xi},
	\end{equation}
with the fewest non-zero terms in $\mathbf{\Xi}$. It is possible to solve for the relevant terms that are active in the dynamics using sparse regression algorithms that penalizes the number of terms in the dynamics and scales well to large problems.  
SINDy works well on low-dimensional systems with a high-quality selection of candidate terms in the library \cite{brunton_discovering_2016, mangan_inferring_2016, champion_data-driven_2019, rudy_data-driven_2017, brunton_sparse_2016}. With too little data, or too many library terms (which can result from too many variables), the algorithm can fail to produce a viable model.

\subsection{Control} 
Control can be broadly divided into the categories of open-loop and closed-loop~\cite{sontag_mathematical_2013, brunton_data-driven_2019}. Most advanced control techniques are closed-loop, based on constructing optimal feedback to stabilize a system at a particular set point.  PID controllers, widely used in industrial control systems, use system error as feedback to efficiently move system output to a desired state without the use of a model of the process \cite{sontag_mathematical_2013, brunton_data-driven_2019}. As PID does not have predictive ability it is not ideal in complex systems with significant latency or higher-order dynamics. Model Predictive Control (MPC), in contrast, uses a dynamical model of the process to predict the effects of an independent variable on the process for a future time window and selects a control signal accordingly \cite{camacho_model_2007, allgower_nonlinear_1999, findeisen_state_2003, kaiser_sparse_2018, brunton_data-driven_2019}. While these methods grant optimal control for systems that can be continuously measured and actuated, many systems of interest are not amenable to constant monitoring or reactive control signals.  Furthermore, highly nonlinear systems may contain many stable attractors that do not require feedback control to stabilize, and can be achieved with transient open-loop control signals \cite{morrison_nonlinear_2020}. We wish to exert control over highly nonlinear systems using open-loop control.  We develop a feedforward control procedure that takes advantage of the nonlinearities of the system to move between stable attractors as well as to create and stabilize fixed points.
\par
An extensive body of literature exists on the analysis and control of nonlinear systems \cite{isidori_nonlinear_1990}.
Linear control techniques can, in many cases, be extended to nonlinear systems using a variety of methods. Controllability criterion, observability, and normal forms have been proposed for nonlinear systems.
Dynamic feedback linearization, sliding mode techniques, and Lyapunov methods can be used to control a variety of nonlinear systems found in chemistry, biology, and electrical and mechanical engineering \cite{isidori_nonlinear_1990}.
Nonlinear model predictive control determines the optimal open-loop control signal trajectories at each sampling instance. Our method is similar to MPC in that we predict the future state of the system under the influence of a control signal using a model.
\par
From a dynamical systems perspective, some research considers the global dynamics that are achievable via feedback control signals by viewing control as parameter manipulations that bring about bifurcations in local dynamics \cite{isidori_nonlinear_1990}. Instead of narrowly considering a single trajectory, this line of research investigates the behavior of ensembles of trajectories. Strategic perturbations to a select number of nodes can induce dynamics on a network to converge to a desired target state, rescuing it from convergence to an undesirable state \cite{cornelius_realistic_2013}. Stochasticity has also been proposed as a useful control mechanism in biological systems. Transcription can be characterized as a bistable dynamical system where stochasticity allowed transcription to fluctuate between the two stable states and therefore acts as a mechanism for regulating transcription \cite{kepler_stochasticity_2001}. Further research found that noise can be exploited to induce desired state transitions in other similarly complex network dynamical systems \cite{wells_control_2015}. Strategic perturbations to tunable parameters is also a form of control as it affects the stability of the nonlinear system’s fixed points.
\par
We continue in this dynamical systems view and consider the global and local dynamics that are possible via system bifurcations induced through feed-forward control signals.
Our control method is limited in that only a subset of all possible states and stability levels can be achieved using feed-forward control in an initial condition agnostic fashion. 
Nonetheless, this type of control may be useful in systems where system measurements are difficult to obtain or feedback mechanisms are impossible.
%
%%%%%%%%%%%%%%%%%%%%%%%%%%%%%%%%%%%%%%%%%%%%%
%%% FF control for low-D systems %%%%%%%%%%%
%%%%%%%%%%%%%%%%%%%%%%%%%%%%%%%%%%%%%%%%%%%%
\section{Feed-forward control for low-dimensional nonlinear systems} 
We consider how to generate control signals that will move a system between attractors (fixed points) in a nonlinear dynamical system $\matr{x'} = F(\matr{x})$, $\matr{x} \in \mathds{R}^2$. We consider systems of the form 
\begin{align}
    \frac{dx}{dt} &= f(x,y) + u_1(t)\\
    \frac{dy}{dt} &= g(x,y) + u_2(t).
\end{align}
where $f$ and $g$ are the intrinsic dynamics and $u_1(t)$ and $u_2(t)$ are the feed-forward control signals. Systems of this form are controlled by regulating the actuating forces.  Many actions taken to control systems can be characterized as feed-forward control signals that do not alter the structure of the underlying system.  Adding reactants to a chemical reaction, removing invasive species, performing deep brain stimulation, and taxing individuals can all be characterized as feed-forward control signals applied to complex dynamical systems.
\par

Biological systems may also utilize feed-forward control signals themselves to modulate internal processes. Behavior transitions in the nematode \textit{C. elegans} can be characterized by feed-forward control signals applied to a nonlinear system with two attractor states \cite{morrison_nonlinear_2020}. While transient feed-forward signals control short-term behavior in this model of \textit{C. elegans}, system parameter modifications affect long-term behaviors, highlighting the different ways biological systems can control their output across multiple timescales. Our control framework is useful both for developing control strategies for nonlinear systems as well as for understanding how natural systems accomplish endogenous control.
%%%%%%%%%%%%%%%%%%%%%%%%%%%%%%%%%%%%%%%%%%%%%5
%%%% stability analysis %%%%%%%%%%%%%%%%%%%%
%%%%%%%%%%%%%%%%%%%%%%%%%%%%%%%%%%%%%%%%%%5%
%
\subsection{Fixed Point Stability}
We use properties of the governing nonlinear system to derive the set of control signals that will move us between stable fixed points by modifying their stability. The Jacobian of the nonlinear system $J(x,y)$ and therefore its trace $T(x,y)$ and determinant $D(x,y)$ are independent of the control signals, meaning that fixed point stability is dependent only on location.  These three quantities are given by the following respectively: 
\begin{align}
    J(x,y) &= \begin{bmatrix}
        f_x(x,y) & f_y(x,y) \\
        g_x(x,y) & g_y(x,y) 
\end{bmatrix} ,\\[0.4em]
 T(x,y) &= f_x(x,y) + g_y(x,y),\\
    D(x,y) &= f_x(x,y) g_y(x,y) - g_x(x,y) f_y(x,y) .
\end{align}
Given a fixed point $(x^*,y^*)$ occurring at $u_1 = -f(x,y)$ and $u_2 = -g(x,y)$. Fixed points exist in one of four stability regions in the trace-determinant plane. Fixed points in the region $T<0$ and $D>0$ are stable sinks, points in the region $T>0$ and $D>0$ are unstable sources, while points in the region $D<0$ are unstable saddles \cite{strogatz_nonlinear_2001}. Note that in the case of a linear system, $J$ is a constant matrix which means that the control signals affect only the location of the fixed point and not its stability. In contrast, feed-forward control signals in a nonlinear system affect both the location and the stability of fixed points; this implies that a much broader range of activity can be generated by feed-forward control signals in nonlinear systems than in linear systems.

Fixed points transition between stability regions along the curves $T(x,y) = 0$ and $D(x,y) = 0$. Transitions between regions can also occur at $D(x,y), \ T(x,y) = \pm \infty$, we therefore define two additional measures: $\hat{D}(x,y) = 1/D(x,y)$ and $\hat{T}(x,y) = 1/T(x,y)$.  We thus note that region transitions also occur along the curves $\hat{D}(x,y) = 0$ and $\hat{T}(x,y) = 0$. Saddle-node bifurcations occur along the curves $D(x,y) = 0$ and $\hat{D}(x,y) = 0$ while Hopf bifurcations occur along the curves $T(x,y) = 0$ and $\hat{T}(x,y) = 0$.

By constraining the relation between $x$ and $y$ by $D(x,y) = 0$ and $\hat{D}(x,y) = 0$ we can solve for parameterized curves for the control signals that eliminate stable fixed points in the system by inducing a saddle-node bifurcation. This curve in the control signal space occurs along
\begin{align}
    C^s(t) = (u_1^s(t), u_2^s(t))
\end{align}
where $u_1^s(t) = -f(t,y(t))$ and $u_2^s(t) = -g(t,y(t))$. The $y(t)$ in these formulae is the implicit solution to $D(t,y(t)) = 0$ or $\hat{D}(t,y(t)) = 0$. In order for $C^s(t)$ to be a boundary between regions, $D$ must switch signs when crossing the curve. $C^s(t)$ is a boundary curve unless $\frac{\partial D}{\partial \bm{v_{\perp}}(t)} = 0$ and $\frac{\partial^2 D}{\partial \bm{v^2_{\perp}}(t)} \neq 0$ where $\bm{v}_{\perp}(t) = -u_2(t) \hat{\matr{i}} + u_1(t) \hat{\matr{j}}$ is the direction orthogonal to $C^s(t)$.

Stable fixed points can also be eliminated through Hopf bifurcations which occur along the curve
\begin{align}
    C^h(t) = (u_1^h(t), u_2^h(t))
\end{align}
where $u_1^h(t) = -f(t,y(t))$, $u_2^h(t) = -g(t,y(t))$, and $y(t)$ is now the implicit solution to $T(t,y(t)) = 0$ or $\hat{T}(t,y(t)) = 0$. $C^h(t)$ is a boundary curve unless $\frac{\partial T}{\partial \bm{v_{\perp}}(t)} = 0$ and $\frac{\partial^2 T}{\partial \bm{v^2_{\perp}}(t)} \neq 0$ where $\bm{v}_{\perp}(t)$ is the direction orthogonal to $C^h(t)$.

\begin{figure}
    \centering
    \includegraphics[width = 1.0\linewidth]{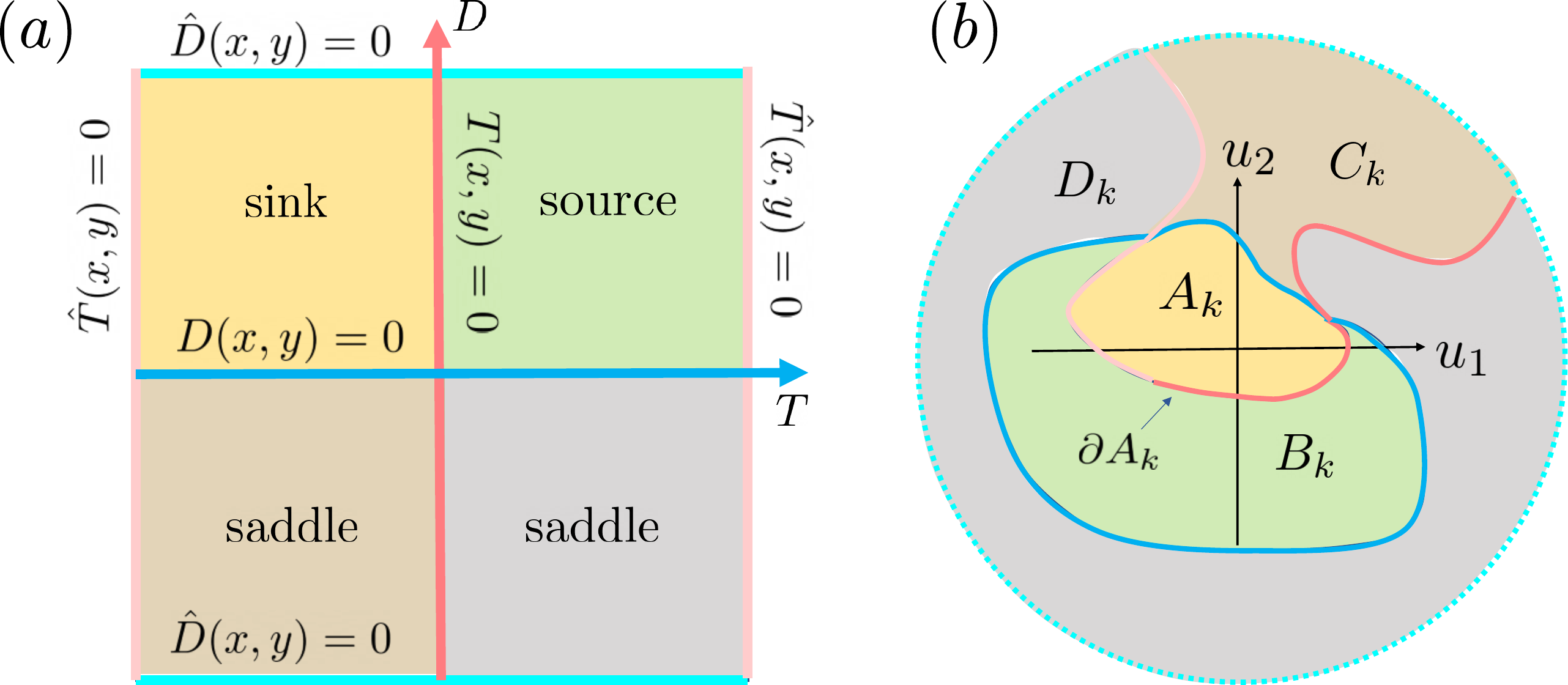}
    \caption{Regions of stability and instability for fixed points in a dynamical system. (a) Fixed points are sinks, sources, or saddles, depending on where they lie in the trace-determinant plane. (b) These regions map to stability regions in the control space. Movement across bifurcation curves in the control space correspond to moving between stability regions in the trace-determinant plane by crossing $D = 0$, $T = 0$, $\hat{D} = 0$, or $\hat{T} = 0$.}
    \label{fig:control_topology}
\end{figure}

We can use these curves to determine stability regions in the control space for each fixed point in the dynamical system.  Let $P = \{p_1, p_2,..., p_n\}$ be the set of fixed points in the nonlinear system $\matr{x'} = F(\matr{x})$. Each fixed point is associated with one of four stability regions $A, B, C,$ or $D$.
Let $A_k$ be the set of control signals $u \in \mathds{R}^2$ such that fixed point $p_k$ is stable under $u$. Let $B_k$ be the set of control signals $u \in \mathds{R}^2$ such that fixed point $p_k$ is a source under $u$. Let $C_k$ be the set of control signals $u \in \mathds{R}^2$ such that fixed point $p_k$ is a saddle with $T<0$ and let $D_k$ be the set of control signals $u \in \mathds{R}^2$ such that fixed point $p_k$ is a saddle with $T>0$.
$A_k$, $B_k$, $C_k$, and $D_k$ are disjoint regions. Fixed point $p_k$ is destabilized in region $A_k^c = \mathds{U} \backslash A_k = B_k \bigcup C_k \bigcup D_k$.
We denote $\partial A_k$ as the boundary of the stability region $A_k$ of fixed point $p_k$. Figure~\ref{fig:control_topology}(a) shows the stability regions and boundaries in the trace-determinant plane while Fig~\ref{fig:control_topology}(b) shows the mapping of these stability regions into the control plane.

%
%
%%%%%%%%%%%%%%%%%%%%%%%%%%%%%%%%%%%%%%%%%%%%%
%%% Controllability given no limit cycles %%%
%%%%%%%%%%%%%%%%%%%%%%%%%%%%%%%%%%%%%%%%%%%%%
\subsection{Controllability Given No Limit Cycles}

We can now define necessary and sufficient conditions for moving between fixed points in a system that does not contain limit cycles. The following conditions concern moving between specific fixed points.
\par
\medskip
\emph{Necessary condition to move directly from $p_\ell$ to $p_k$:}
\begin{align}
    A_{\ell}^c \bigcap A_k \neq \emptyset
\end{align}
\par
\medskip
\emph{Sufficient conditions to move directly from $p_\ell$ to $p_k$:}
\begin{align}
    \bigcap_{i=1, i \neq k}^n A_i^c \bigcap A_k \neq \emptyset
\end{align}
\par
We can extend these conditions to make statements about the reachability of all fixed points in the system.
\par
\medskip
\emph{Necessary condition to move directly from any fixed point in the system to any other fixed point:}
\begin{align}
    \forall \ell,k \in 1,...,n \quad A_{\ell}^c \bigcap A_k \neq \emptyset
\end{align}
\par
\medskip
\emph{Sufficient conditions to move directly from any fixed point in the system to any other fixed point:}
\begin{align}
    \forall k \in 1,...,n  \bigcap_{i=1, i \neq k}^n A_i^c \bigcap A_k \neq \emptyset
\end{align}
\par
\medskip
\emph{Necessary conditions to move to any point in the system:}
\begin{align}
    \forall k \in 1,...,n \ \exists \ell \in 1,...,n \ s.t. \ A_\ell^c \bigcap A_k \neq \emptyset
\end{align}
We can use the control regions to move between fixed points in the system by destabilizing fixed points we want to escape, and stabilizing fixed points we want to achieve.
Transitioning between fixed points through saddle-node bifurcations is preferable to transitioning via Hopf bifurcations as Hopf bifurcations can create a stable limit cycle around the source.
\par
While these sets generate collections of control signals that can be used to move between fixed points, it does not reveal the optimal control signals for a transition which depends on additional goals such as energy minimization, transition speed, robustness, or preferred path.

%%%%%%%%%%%%%%%%%%%%%%%%%%%%%%%%%%%%%%%
%%%%% Limit Cycles %%%%%%%%%%%%%%%%%%%%
%%%%%%%%%%%%%%%%%%%%%%%%%%%%%%%%%%%%%%%
\subsection{Limit Cycles}
\label{sec:limit_cycles}
\par
Limit cycles are a pervasive feature of many nonlinear systems and so must be considered.   Unfortunately their location and stability are more difficult to characterize than the location and stability of fixed points in the system. We can make some statements about where limit cycles do and do not appear by using Dulac's criterion and the Poincar\'e-Bendixson theorem \cite{strogatz_nonlinear_2001}.
We can use Dulac's criterion to solve for regions encapsulating fixed points that cannot contain limit cycles and we can use the Poincare-Bendixson theorem to solve for regions that do contain limit cycles.
For example, a globally stable system with a single fixed point that is a source is a compact set and therefore must contain at least one stable limit cycle encapsulating the fixed point.
Some systems may contain compact sets enclosing all fixed points. Such sets can be found by eliminating the highest order terms that determine the global stability and then observing the stability of the resulting system as $t \rightarrow \infty$.  Thus we consider
\begin{align}
    x' \approx P_{n-1}(x,y) + p(x^n, y^n)\\
    y' \approx Q_{n-1}(x,y) + q(x^n, y^n)
\end{align}
If the global stability of $(P_{n-1}, Q_{n-1})$ is opposite the global stability of $(f, g)$ and along a Jordan curve surrounding all fixed points $|P_{n-1}(x,y)|> |p(x^n,y^n)|$ and $|Q_{n-1}(x,y)|> |q(x^n,y^n)|$, then there must be a compact set outside of this Jordan curve and therefore a limit cycle.

In some systems multiple limit cycles may encapsulate a fixed point or region. We may find the presence of multiple layers of limit cycles by repeating the process on the approximate lower-order system $P_{n-1}(x,y) = P_{n-2}(x,y) + p(x^{n-1},y^{n-1})$ and $Q_{n-1}(x,y) = Q_{n-2}(x,y) + q(x^{n-1},y^{n-1})$. Once again if the lower order terms are larger in absolute value than the higher order
terms along a Jordan curve surrounding the fixed points and enclosed by the first Jordan curve, then there must be a second limit cycle surrounding the fixed points.
%%%%%%%%%%%%%%%%%%%%%%%%%%%%%%%%

We can determine the global stability of the system $\matr{x}' = F(\matr{x})$, that is, the stability as $x,y \rightarrow \pm \infty$, by instead considering the stability of the system $\hat{\matr{x}}' = G(\hat{\matr{x}})$ as $\hat{\matr{x}} \rightarrow 0$ where $\hat{x} = 1/x$ and $\hat{y} = 1/y$,

\begin{align}
    \frac{d\hat{x}}{dt} &= -\hat{x}^2 \left( f \left( \frac{1}{\hat{x}},\frac{1}{\hat{y}} \right) + u_1(t) \right)\\
    \frac{d\hat{y}}{dt} &= -\hat{y}^2 \left( g \left(\frac{1}{\hat{x}},\frac{1}{\hat{y}} \right) + u_2(t) \right)
\end{align}

Our original system $\matr{x}' = F(\matr{x})$ is stable as $x,y \rightarrow \pm \infty$ if $G(\hat{\matr{x}})$ does not have any stable directions as $\hat{x}, \hat{y} \rightarrow 0$ and unstable otherwise.

%%%%%%%%%%%%%%%%%%%%%%%%%%%%%%%%%5

We can determine that there are no limit cycles surrounding all fixed points by mapping the system to $G(\hat{\matr{x}})$ and using the Dulac's criterion to find a region surrounding $\hat{\matr{x}}=0$ that does not contain any limit cycles.
Conversely, we can determine if there are limit cycles surrounding a particular fixed point in a system with multiple fixed points by mapping that point to infinity in the system $G(\hat{\matr{x}})$ and then finding if there is a limit cycle surrounding all fixed points in $G(\hat{\matr{x}})$ using the method outlined above with the Poincare-Bendixson theorem.

\subsection{Controllability with Limit Cycles}
We can use bifurcations that create and eliminate limit cycles to move between limit cycles and fixed points in the system. Andronov-Hopf bifurcations eliminate stable limit cycles by turning the source in the center of the limit cycle into a sink \cite{kuznetsov_elements_1998, guckenheimer_nonlinear_1983}. Homoclinic saddle-node bifurcations, otherwise known as infinite period bifurcations, create a saddle-node bifurcation along the limit cycle which eliminates the cycle
\cite{kuznetsov_elements_1998,  guckenheimer_nonlinear_1983, keener_infinite_1981}. Homoclinic bifurcations occurs when a limit cycle merges with a saddle point creating a homoclinic orbit.  Saddle-node bifurcations of periodic orbits merge concentric stable and unstable limit cycles \cite{strogatz_nonlinear_2001}.
\par
Feed-forward control may be able to create some of these bifurcations in a given system. We can control the creation and elimination of limit cycles in the dynamical system by using Dulac's criterion and the Poincare-Bendixson theorem to map out control regions that will create and eliminate limit cycles.
\subsection{Extension to Three-dimensional Systems}
This analysis can be extended to create control procedures for three-dimensional systems, $\matr{x}' = F(\matr{x}) + \matr{u}(t)$, where $\matr{x}, \matr{u} \in \mathds{R}^3$. 
\begin{align}
    \frac{dx}{dt} = f(x,y,z) + u_1(t)\\
    \frac{dy}{dt} = g(x,y,z) + u_2(t)\\
    \frac{dz}{dt} = h(x,y,z) + u_3(t) .
\end{align}
Analogous to the two-dimensional case, we compute the Jacobian of the system $J(x,y,z)$ and fixed point locations $u_1 = -f(x,y,z)$, $u_2 = -g(x,y,z)$, $u_1 = -h(x,y,z)$. We use the linearized system in conjunction with the control signal equations to determine the control surface $(u_1(t,s), u_2(t,s), u_3(t,s))$ along which bifurcations occur. While the stability of fixed points can still be determined in three-dimensional systems, the locations of strange attractors become impossible to determine analytically as
Dulac's criterion and the Poincare-Bendixson theorem only apply to two-dimensional systems \cite{strogatz_nonlinear_2001}. Therefore, while we can still stabilize and destabilize fixed points in order to execute transitions, it is unknown where and when attractors will occur in the system, thus potentially compromising the control procedure advocated here. In these circumstances, attractors and their bifurcations are best found experimentally by simulating the dynamics. 
Using data, Poincar\'e maps can be found via the SINDy algorithm and from these discovered maps the stability of nonlinear periodic orbits can be determined \cite{bramburger_poincare_2019}.

%%%%%%%%%%%%%%%%%%%%%%%%%%%%%%%%%%%%%%%%%%%%%%%%%%
%%% Section - low-D control applications %%%%%%%%%
%%%%%%%%%%%%%%%%%%%%%%%%%%%%%%%%%%%%%%%%%%%%%%%%%%
\section{Network control applications}
Our feed-forward control procedure can be directly applied to nonlinear systems that consist of two or three variables. While most real world systems realistically involve a large number of variables, many canonical models in chemistry, biology, and the social sciences are formulated using only a few variables.
In chemistry the Brusselator describes an autocatalytic chemical reaction with two reagents \cite{prigogine_symmetrybreaking_1967, prigogine_symmetry_1968, strogatz_nonlinear_2001} and the Lokta-Volterra model in ecology describes predator-prey interactions at the population level between two species \cite{lotka_analytical_1920, lotka_elements_1925, volterra_fluctuations_1926}.
In epidemiology the SIR model describes the spread of an infectious disease by modeling the interactions between susceptible and infected individuals \cite{kermack_contribution_1927, hethcote_mathematics_2000} and in neuroscience the FitzHugh–Nagumo model describes spike generation in neurons by measuring membrane voltage and a recovery variable \cite{fitzhugh_impulses_1961, nagumo_active_1962}.
In the social sciences the classical model of political economy uses capital and labor as variables \cite{samuelson_canonical_1978} while Richardson's arms race model describes the interactions between two nations stockpiling nuclear weapons \cite{richardson_arms_1960, smith_influence_2020}.
If the dynamics of a low-dimensional system are unknown, they can be found using the SINDy algorithm \cite{brunton_discovering_2016, brunton_data-driven_2019}, allowing us to apply control even to novel systems.
We demonstrate the effects of feed-forward control on two variable nonlinear systems using simple chemical reaction models as examples.

%%%%%%%%%%%%%%%%%%%%%%%%%%%%%%%%%%%%%%%%%%%%%
%% Subsection - Chem reaction %%%%%%%%%%%%%%%
%%%%%%%%%%%%%%%%%%%%%%%%%%%%%%%%%%%%%%%%%%%%%
\subsection{Chemical Reaction with Bistability} 
\begin{figure}[t]
    \centering
    \includegraphics[width=0.9\linewidth]{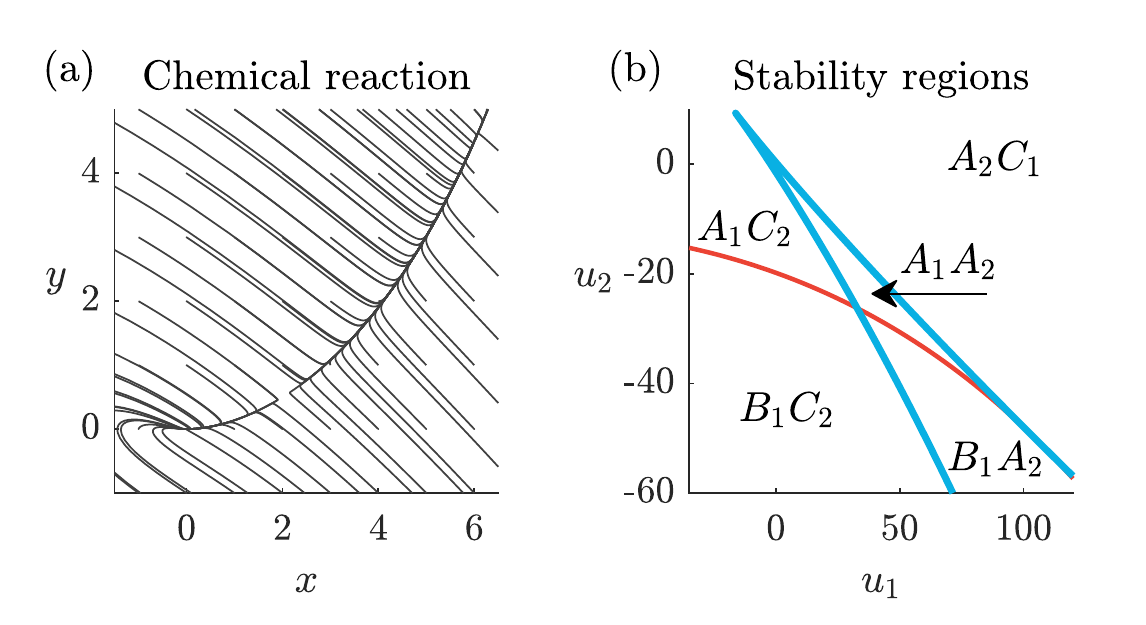}
    \caption{Chemical reaction with bistability. (a) In the absence of control signals there are two stable fixed points. (b) Regions of stability and instability for the fixed points in the uncontrolled system. The system has regions with two stable fixed points, one stable fixed point, a stable fixed point and a stable limit cycle, and only a stable limit cycle. Note the stable limit cycle appears surrounding the unstable source.}
    \label{fig:chem_reaction}
\end{figure}

The following is a minimal example of a chemical reaction with bistability using parameter values $k_1 = 8$, $k_2 = 1$, $k_3 = 1$, and  $k_4 = 3/2$ from \cite{wilhelm_smallest_2009}. We can find control signals that will move the system between stable states by computing stability regions from the control signal bifurcation curves.  The governing equations are given by
\begin{align}
    \frac{dx}{dt} &= 16y - x^2 - xy -\frac{3}{2}x + u_1\\
    \frac{dy}{dt} &= x^2 - 8y + u_2
\end{align}
where the Jacobian, trace and determinant are given by
\begin{align}
    J(x,y) = \begin{bmatrix}
        -2x - y - \frac{3}{2} & 16-x \\
        2x & -8
\end{bmatrix}
\end{align}
\begin{align}
    T(x,y) &= -2x  - y -\frac{19}{2}\\
    D(x,y) &= 2x^2 - 16x +8y + 12 .
\end{align}
The curve for $T(x,y) = 0$ is
\begin{align}
    u_1^h(t) &= -t^2+24t+152\\
    u_2^h(t) &= -t^2 - 16t-76
\end{align}
The curve for $D(x,y) = 0$ is
\begin{align}
    u_1^s(t) &= -\frac{1}{4}t^3 + 7t^2 - 32t + 24\\
    u_2^s(t) &= -3t^2 + 16t - 12
\end{align}
In the uncontrolled system the chemical reaction has two stable fixed points, Fig.~\ref{fig:chem_reaction}(a). Figure~\ref{fig:chem_reaction}(b) shows the stability regions of the fixed points in the control space. Across the blue curve fixed points exhibit saddle-node bifurcations while across the red curve fixed points exhibit Hopf bifurcations. This map shows us that we can eliminate the second fixed point via a saddle-node bifurcation while keeping the first fixed point stabilized and then destabilize the first fixed point via a Hopf bifurcation which creates a stable limit cycle surrounding the first fixed point.
The overlay of stability regions shows us the wide variety of states that can be achieved by simply adding or removing certain quantities of reactant.

%%%%%%%%%%%%%%%%%%%%%%%%%%%%%%%%%%%%%%%%%%
%%% Subsection - Brusselator %%%%%%%%%%%%%
%%%%%%%%%%%%%%%%%%%%%%%%%%%%%%%%%%%%%%%%%%
\subsection{Brusselator}
The Brusselator describes a type of autocatalytic chemical reaction. We take $a = 1$ and $b = 3$ from the general model \cite{strogatz_nonlinear_2001} 
\begin{align}
    \frac{dx}{dt} &= 1 + x^2 y - 4x+ u_1(t)\\
    \frac{dy}{dt} &= 3x - x^2 y + u_2(t) ,
\end{align}
and find the bifurcation curve locations.

The Jacobian of the system is
\begin{align}
    J(x,y) = \begin{bmatrix}
        2xy-4 & x^2\\
        3-2xy & -x^2
\end{bmatrix}
\end{align}
with trace and determinant given by
\begin{align}
    T(x,y) &= 2xy - x^2 -4\\
    D(x,y) &= x^2.
\end{align}
We find that 
$T(x,y) = 0$ along the curve
\begin{align}
    u_1^h(t) &= \frac{-t^3}{2}+2t - 1\\
    u_2^h(t) &= \frac{t^3}{2}-t .
\end{align}
Note that $\hat{T}(x,y) \rightarrow \pm 0$ when $y = c/x^2$ and $x \rightarrow \pm 0$ which produces the curve
\begin{align}
    u_1^h(t) &= -1 -t\\
    u_2^h(t) &= t
\end{align}
The determinant is always greater than zero and therefore does not produce a curve over which there is a sign change. This means that feed-forward control signals that induce saddle-node bifurcations do not exist in this system.
\par
\begin{figure}[t]
    \centering
    \includegraphics[width=1.0\linewidth]{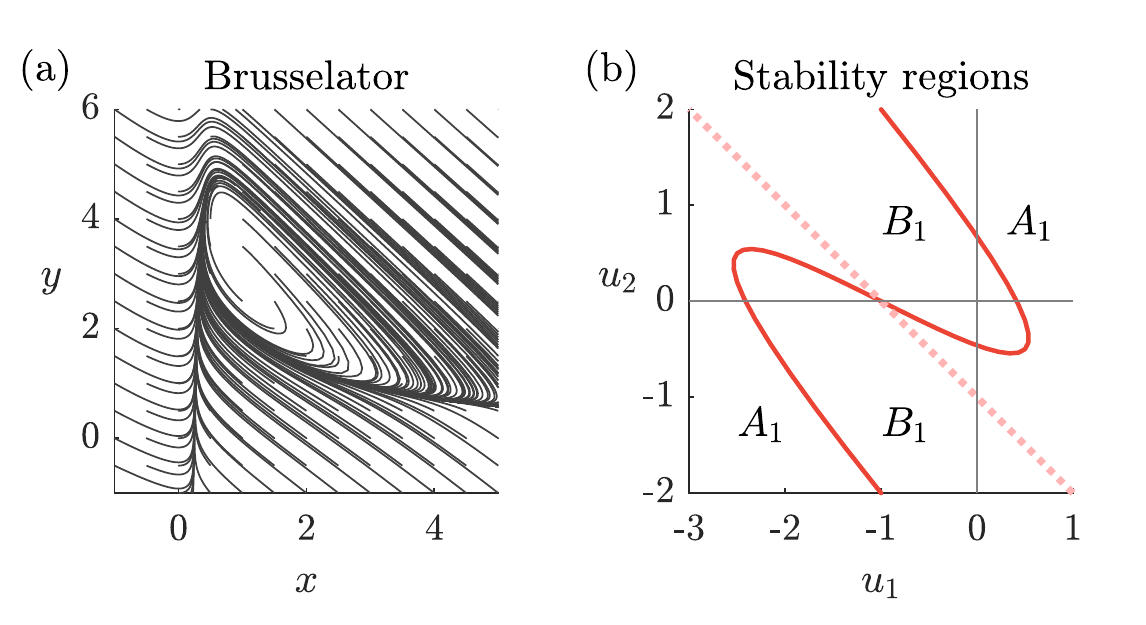}
    \caption{Brusselator (a) Brusselator with $u_1, u_2 = 0$. (b) Stability regions in the control space. The system transitions between having a single stable fixed point and an unstable source encapsulated by a stable limit cycle.}
    \label{fig:brusselator}
\end{figure}

Figure~\ref{fig:brusselator}(a) shows the uncontrolled Brusselator contains a stable limit cycle as also indicated by the control signal stability regions, Fig.~\ref{fig:brusselator}(b). Across the Hopf bifurcation curve $C^h(t)$ the globally stable system transitions between having a single stable fixed point and a single source surrounded by a limit cycle, Fig.~\ref{fig:brusselator}(b). The Brusselator can be compelled to move between a stable state and a limit cycle by adding or removing particular amounts of reactant.
%
%%%%%%%%%%%%%%%%%%%%%%%%%%%%%%%%%
%%% Section - High-D control %%%%
%%%%%%%%%%%%%%%%%%%%%%%%%%%%%%%%%
\section{Feed-forward control for high-dimensional nonlinear systems}
So far we have considered control for systems with only two or three variables. While some systems fit this description, most real systems involve a much larger number of variables and interactions.
We can use dimensionality reduction in conjunction with data-driven discovery of nonlinear dynamics to extend this control technique to high-dimensional systems.
Consider an n-dimensional nonlinear system of the form 
\begin{align}
    \matr{x}' = G(\matr{x}) + \matr{u}_G(t)
\end{align}
where the system variables $\matr{x} \in \mathds{R}^n$ are controlled by the feed-forward control signal $\matr{u}_G \in \mathds{R}^n$ and the form of $G$ may be unknown. We detail a method for finding $\matr{u}_G(t)$ from system measurement data, as outlined in Algorithm~\ref{alg:highD}.  The steps of the algorithm are detailed as follows:
%
%%%%%%%%%%%%%%%%%
%%%%%%%%%%%%%%%%%%%%%
\begin{figure}[!t]
 \removelatexerror
  \begin{algorithm}[H]
    \label{alg:highD}
   \caption{Control for high dimensional systems}
   Measure system outputs: $\matr{X}$\\
   Dimensionality reduction: $\matr{U} \matr{S} \matr{V^*}$\\
   Identify nonlinear system: $\matr{x'} = F(\matr{x})$\\
   Determine control signals: $\matr{x'} = F(\matr{x}) + \matr{u}_{F}(t)$\\
   Determine high dimensional signals: $\matr{u}_{G}(t) = \matr{\hat{U}}\matr{u}_{F}(t)$
  \end{algorithm}
\end{figure}
%%%%%%%%%%%%%%5
%%%%%%%%%%%%%%%%
%
\subsubsection{Measure system outputs} 
We begin by collecting timeseries data from variables in the high-dimensional system. Data should be collected from many different initial conditions in the entire state space in order to generate an accurate model.

\subsubsection{Dimensionality reduction}
Next we dimensionality reduce the timeseries data using the SVD. Not all high-dimensional systems have a low-dimensional structure and those that do have a low-dimensional structure may not necessarily have a \emph{linear} low-dimensional structure. SVD is a suitable dimension reduction technique for a given dataset if the first two modes capture the majority of the variance in the system.

\subsubsection{Identify the nonlinear system using SINDy}
Once we have reduced the data to only a few dimensions we can fit a nonlinear dynamical system to the data using the SINDy algorithm \cite{brunton_discovering_2016}. If the original system is globally stable the model should also be globally stable. We can determine the stability of the SINDy model using the criteria outlined in Section~\ref{sec:limit_cycles}.

\subsubsection{Determine feed-forward control signals}
Now that we have a low-dimensional nonlinear model for the system's activity, we can generate feed-forward control signals that stabilize fixed points and transition the system between attractors in the SINDy model. These control signals must go through a last transformation before they can be applied to the original system.

\subsubsection{Determine high dimensional control signals}
Finally, we use the feed-forward control signals $\matr{u}_F(t)$ discovered for the low-dimensional system to control the original system by projecting the signals back to the original high-dimensional space via the first two SVD modes $\matr{\hat{U}}$.
\begin{align}
    \matr{u}_G(t) &= \matr{\hat{U}}\matr{u}_F(t)\\
    \matr{x'} &= G(\matr{x}) + \matr{u}_G(t) .
\end{align}
Developing control methods in a reduced space has been effective in other settings \cite{brunton_data-driven_2019}. Dynamic mode decomposition with control (DMDc) finds spatial-temporal coherent modes with which to construct low-order models that incorporate the effects of control signals \cite{proctor_dynamic_2016, kutz_dynamic_2016}. Unsteady wake flows can be effectively described and controlled using POD modes \cite{nair_networked-oscillator-based_2018}. Structural balance dynamics in social networks and the bifurcations that occur in the system can be analyzed in a low-dimensional eigenspace \cite{marvel_continuous-time_2011, morrison_community_2019}.
Similar to previous methods, we use linear dimension reduction, but unlike previous methods we employ nonlinear, feed-forward control.
%
%%%%%%%%%%%%%%%%%%%%%%%%%%%%%%%%%%%%%%%%%%%%%%%%%%
%%% Section - High-D control applications %%%%%%%%
%%%%%%%%%%%%%%%%%%%%%%%%%%%%%%%%%%%%%%%%%%%%%%%%%%
\section{High-dimensional control applications}
High-dimensional systems that exhibit low-dimensional dynamics are pervasive in nature. The nematode \textit{C. elegans} has a network of 302 neurons yet exhibits neural activity that exists on a three-dimensional manifold \cite{kato_global_2015,fieseler_unsupervised_2020}. The low-dimensional encoding of information is a common motif found throughout the neuroscience literature \cite{pfau_robust_2013, briggman_optical_2005, machens_functional_2010, stopfer_intensity_2003} as well as in the study of artificial neural networks \cite{farrell_recurrent_2019, recanatesi_dimensionality_2019, hinton_reducing_2006, jaderberg_speeding_2014, denton_exploiting_2014, rigamonti_learning_2013}. In the social sciences, social networks and structural balance dynamics both exhibit low-dimensional structures \cite{belkin_laplacian_2002, tang_line_2015, marvel_continuous-time_2011, morrison_community_2019}. In biology and epidemiology the dimension of population models can be reduced through mathematical aggregation methods \cite{auger_aggregation_2008}. In molecular biology the dynamics of biological networks exhibiting multiple timescales can be simplified by modeling the system with multiple timescale networks \cite{wang_modelling_2004} and reduced, hierarchical structured models can describe the dynamics of complex signal transduction networks \cite{conzelmann_domain-oriented_2006}. 
We demonstrate our control procedure on several nonlinear high-dimensional systems that exhibit dynamics on a low-dimensional manifolds.

%%%%%%%%%%%%%%%%%%%%%%%%%%%%%%%%%%%%%%%%
%%% Subsection - Random networks %%%%%%%
%%%%%%%%%%%%%%%%%%%%%%%%%%%%%%%%%%%%%%%%
\subsection{Random High-dimensional Networks}
\begin{figure}[t]
    \centering
    \includegraphics[width=0.95\linewidth]{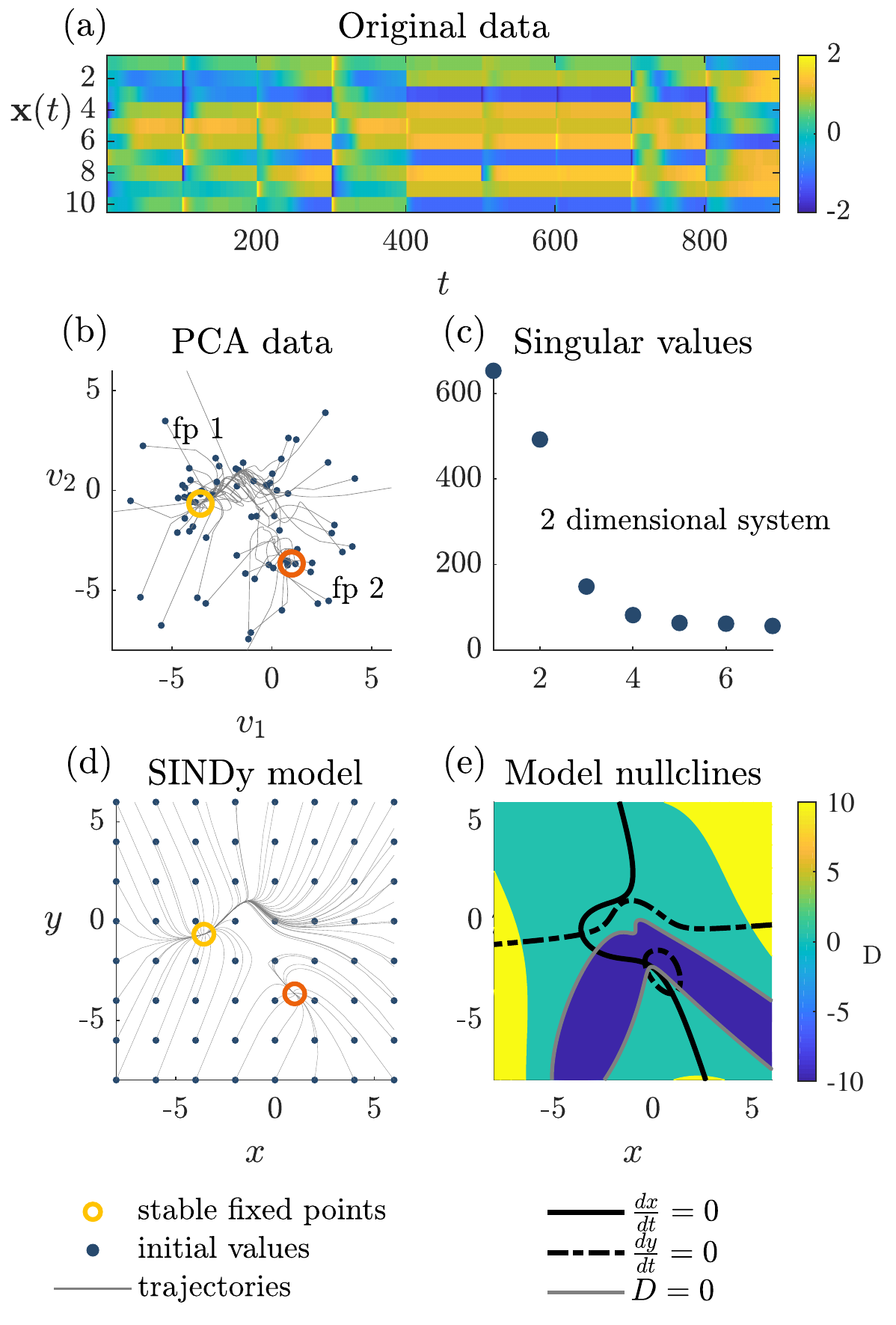}
    \caption{Low-dimensional model of a high dimensional random dynamical system.
    (a) High dimensional system $n=10$ visualized as a timeseries and (b) shown in PCA space. 
    (c) The first two PCA modes capture a reasonable amount of the data and the SINDy model (d) finds the two stable fixed points.
    (e) Nullclines and determinant of the SINDy model.}
    \label{fig:random2}
\end{figure}

We first illustrate our feed-forward control technique on a randomly generated high-dimensional network dynamical system. We wish to find stable attractors endogenous to the system and then generate control signals that will move the system between stable states.
Figure~\ref{fig:random2}(a) shows data collected from our high-dimensional random network. The trajectories in the underlying low-rank SVD/PCA space indicates that there are two stable fixed points in the system. Figure~\ref{fig:random2}(b) shows the system is indeed low-dimensional as indicated by the rapid decay in singular values in Fig.~\ref{fig:random2}(c). The SINDy algorithm finds a model that represents the low-dimensional system and which can be used to derive control signals Fig.~\ref{fig:random2}(d-e).
\par
\begin{figure}[t]
    \centering
    \includegraphics[width=1.05\linewidth]{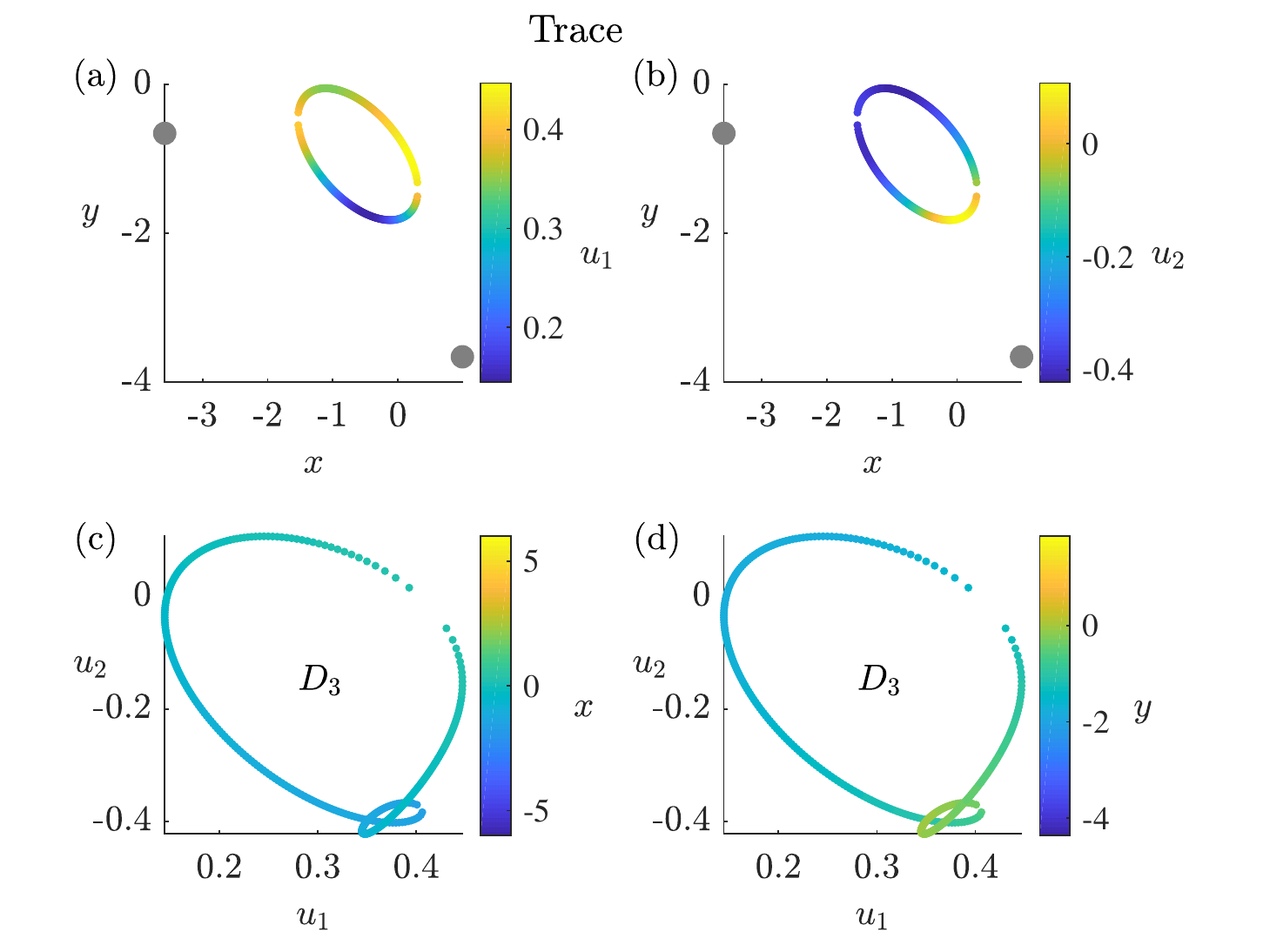}
    \caption{Fixed point locations and control signal values along the curve $C^h(t)$ for the low-dimensional model found in Figure~\ref{fig:random2}.
    (a-b) Fixed point locations in the xy-plane for $C^h(t)$ colored by the $u_1, u_2$ control signal values along $C^h(t)$.
    (c-d) Control signal locations in the $u_1,u_2$ plane for $C^h(t)$ colored by $(x,y)$ fixed point location values.
    The fixed point in this model that transitions across the curve $T = 0$ is the saddle fixed point that sits between the stable fixed points in the uncontrolled model. This saddle fixed point originally is located in region $C_3$ but transitions to region $D_3$ across the curve $C^h(t)$.}
    \label{fig:random2_bifs_trace}
\end{figure}

We next use the SINDy model to find the control signal curves along which bifurcations in the system occur, $C^h(t)$ and $C^s(t)$. Figure~\ref{fig:random2_bifs_trace} shows the locations of fixed points and control signals in the system when $T(x,y) = 0$. Figure~\ref{fig:random2_bifs_trace}(a-b) shows the location of a saddle fixed point as it switches types colored by the control signal values that induce these bifurcations. Figure~\ref{fig:random2_bifs_trace}(c-d) shows these control signal values colored by fixed point locations. The control signals move from stability region $C$ to $D$ as they cross the bifurcation curve.
\par
\begin{figure}[t]
    \centering
    \includegraphics[width=1.05\linewidth]{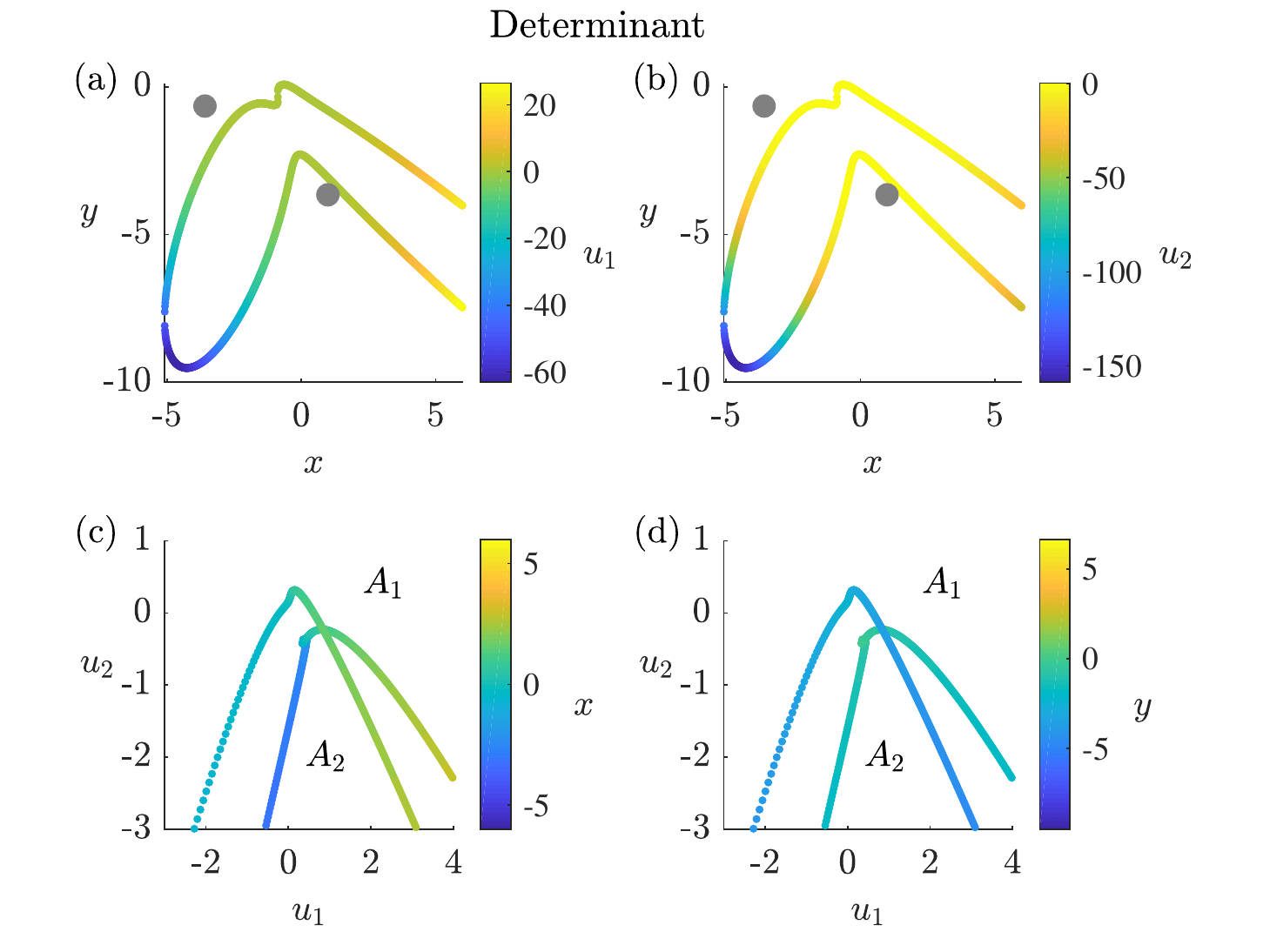}
    \caption{Fixed point locations and control signal values along the curve $C^s(t)$ for the low-dimensional model found in Figure~\ref{fig:random2}.
    (a-b) Fixed point locations in the xy-plane for $C^s(t)$ colored by the $u_1, u_2$ control signal values along $C^s(t)$. 
    (c-d) Control signal locations in the $u_1, u_2$ plane for $C^s(t)$ colored by $(x,y)$ fixed point location values. The fixed points in this model that transition across the curve $D = 0$ are the stable fixed points in the uncontrolled model. The right stable fixed point goes through a saddle-node bifurcation along the top curve, while the left stable fixed point goes through a saddle-node bifurcation along the bottom curve.}
    \label{fig:random2_bifs_det}
\end{figure}

Figure~\ref{fig:random2_bifs_det} shows the locations of fixed points and control signals in the system when $D(x,y) = 0$. Figure~\ref{fig:random2_bifs_det}(a-b) show the locations of the stable fixed points as they go through saddle-node bifurcations colored by the control signal values that induce these bifurcations while Figure~\ref{fig:random2_bifs_det}(c-d) show these control signal values colored by fixed point locations. The stable fixed points start in stability region $A$ and then disappear as they cross the bifurcation curves. $C^s(t)$ defines stability region borders which give us the ability to select control signals that will induce transitions between stable states.

\par
\begin{figure}[!ht]
    \centering
    \includegraphics[width=1\linewidth]{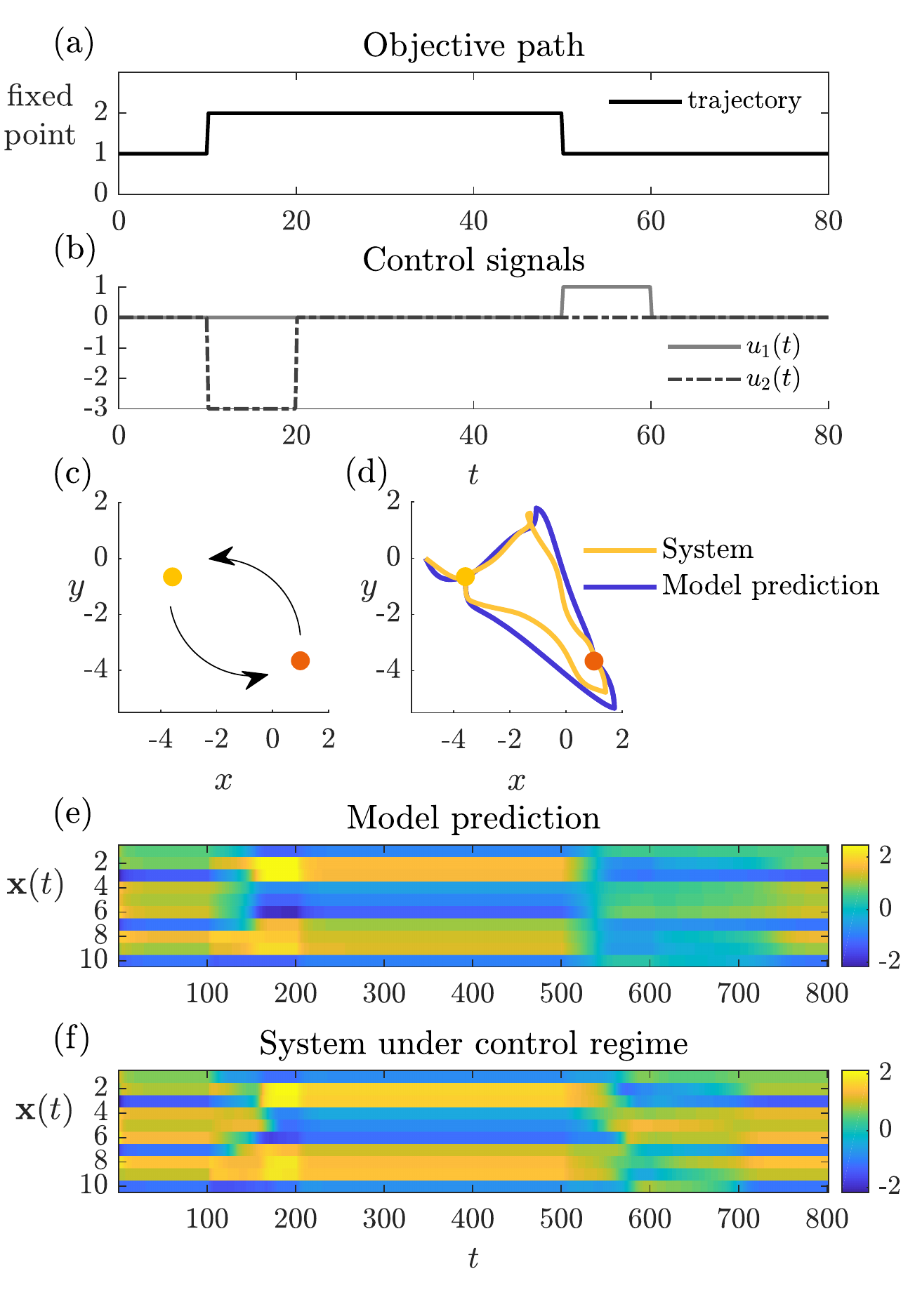}
    \caption{Control example for the random system in Figure~\ref{fig:random2}.
    (a) Objective path for the system.
    (b) Control signals selected using the system's stability region maps to move the system between fixed points.
    (c) Objective path in PCA space.
    (d) Predicted and actualized system paths in PCA space.
    (e-f) Predicted and actualized system activity in the original high-dimensional space.}
    \label{fig:random2_controlled}
\end{figure}

We can now induce transitions between stable states in the system using transient feed-forward control signals selected from the generated stability regions. We first stipulate our objective path or what can be thought of as setpoints over time Fig.~\ref{fig:random2_controlled}(a) and then select control signals that will induce each transition Fig.~\ref{fig:random2_controlled}(b). When the control signals are applied to the original system it makes the desired transitions as viewed in the PCA space Fig.~\ref{fig:random2_controlled}(c-d) as well as the high-dimensional space Fig.~\ref{fig:random2_controlled}(e-f). The system does not behave exactly as the SINDy model predicts as the model is only an approximation of the system's true dynamics.
%

%%%%%%%%%%%%%%%%%%%%%%%%%%%%%%%%%%%%%%%
%%% Subsection - Hopfield %%%%%%%%%%%%%
%%%%%%%%%%%%%%%%%%%%%%%%%%%%%%%%%%%%%%%

\subsection{Hopfield Networks}
Hopfield networks are auto-associative memory networks that converge to stable patterns in the presence of noise and are a model of memory retrieval in the human brain \cite{hopfield_neural_1982, morrison_preventing_2017}. The heaviside function is a key component of the Hopfield model; it is the nonlinearity in the system that allows the Hopfield network to have many fixed points as opposed to one. Because the Heaviside function is a dimension reducing function it makes the dynamics of the Hopfield model much lower dimensional than the network's original size, making it a promising candidate for our control procedure.

\begin{figure}[t]
    \centering
    \includegraphics[width=0.95\linewidth]{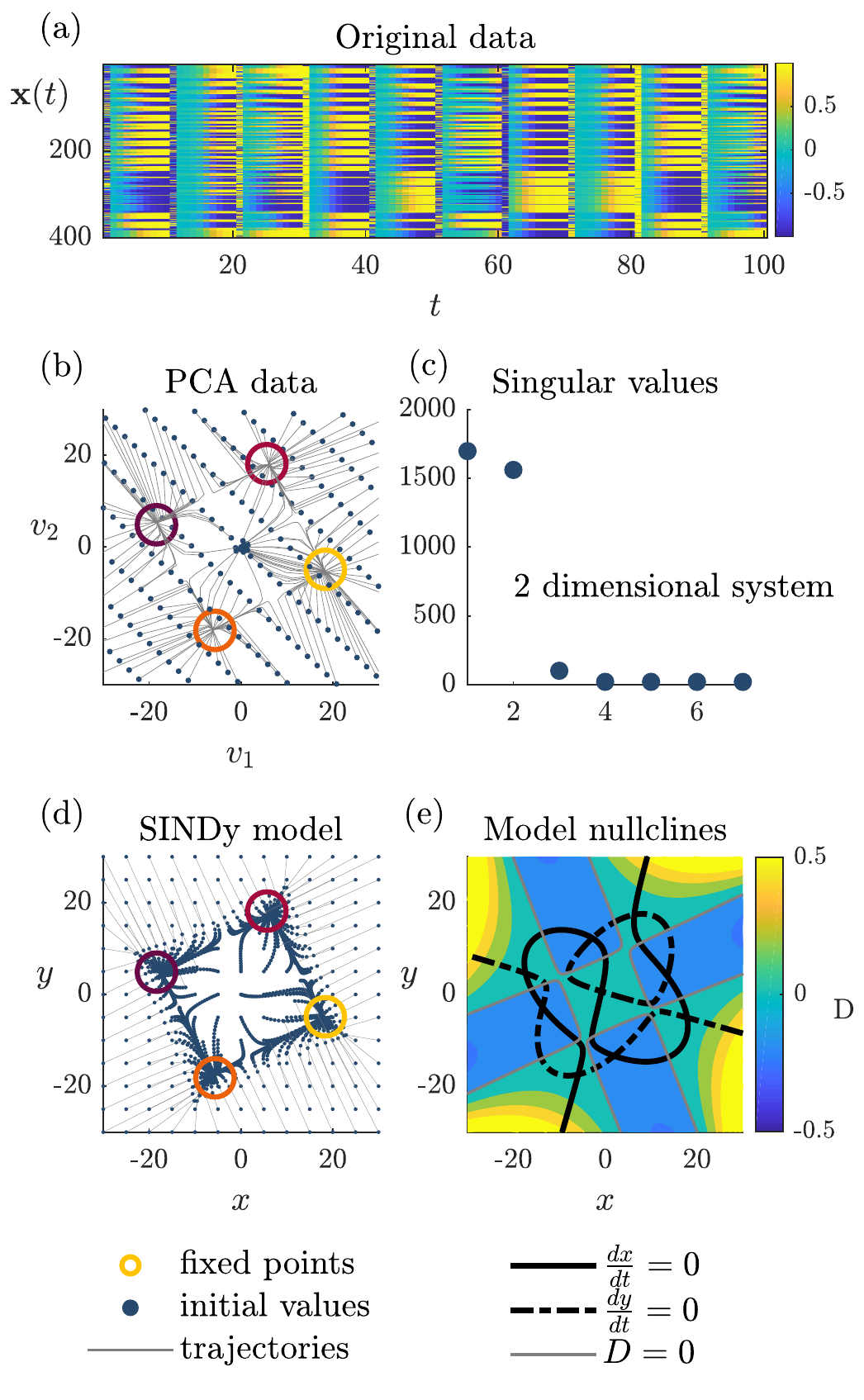}
    \caption{Smooth Hopfield network size $n=400$ with intrinsic 2-dimensional structure. (a) High-dimensional data measured from many initial conditions.
    (b) Data in PCA space. The network converges to one of 4 "memories".
    (c) The first 2 modes capture the majority of the variance in the system. (d) SINDy is able to generate a low-dimensional model with similar dynamics that we can use to generate stability regions.
    (e) Nullclines and determinant of the SINDy model.}
    \label{fig:hop2}
\end{figure}

\begin{figure}[t]
    \centering
    \includegraphics[width=0.7\linewidth]{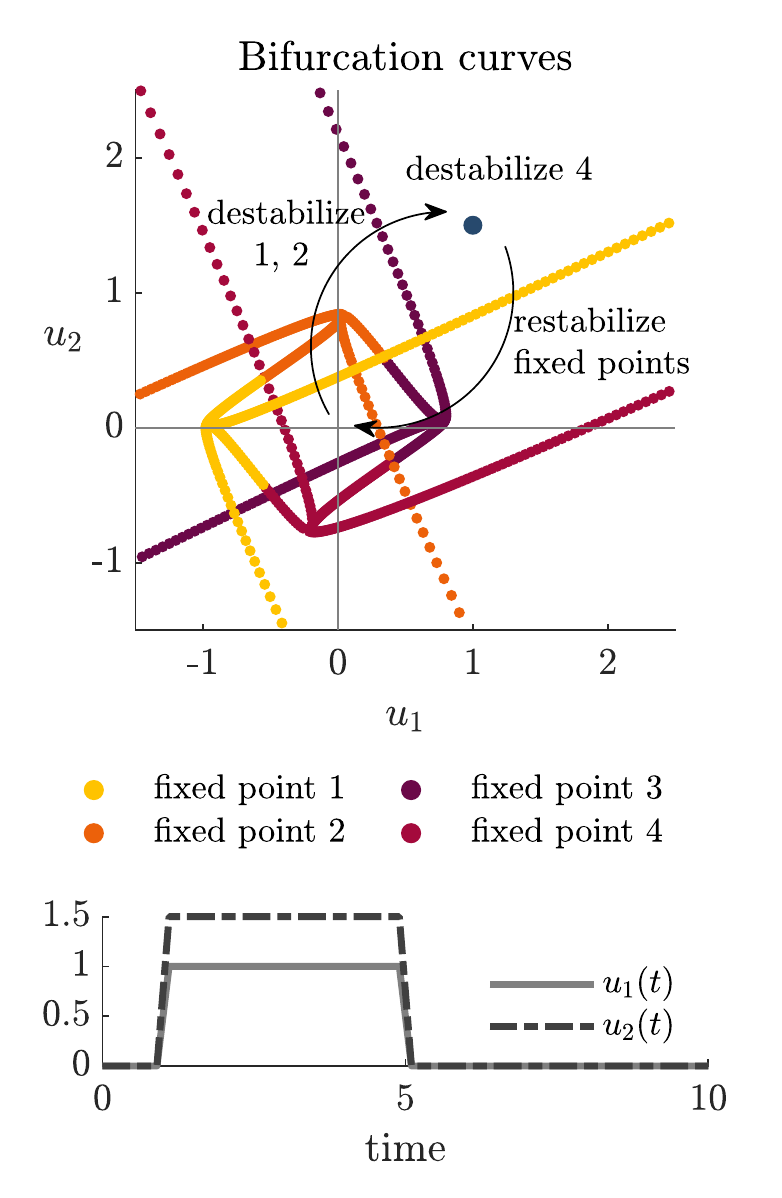}
    \caption{Saddle-node bifurcation curves $C^s(t)$ colored by fixed point for the system in Figure~\ref{fig:hop2}. Moving across the saddle-node bifurcation curve of a stable fixed point either eliminates it from the system or reinstates it.}
    \label{fig:hopfield2_1_bif}
\end{figure}

\begin{figure}[t]
    \centering
    \includegraphics[width=0.9\linewidth]{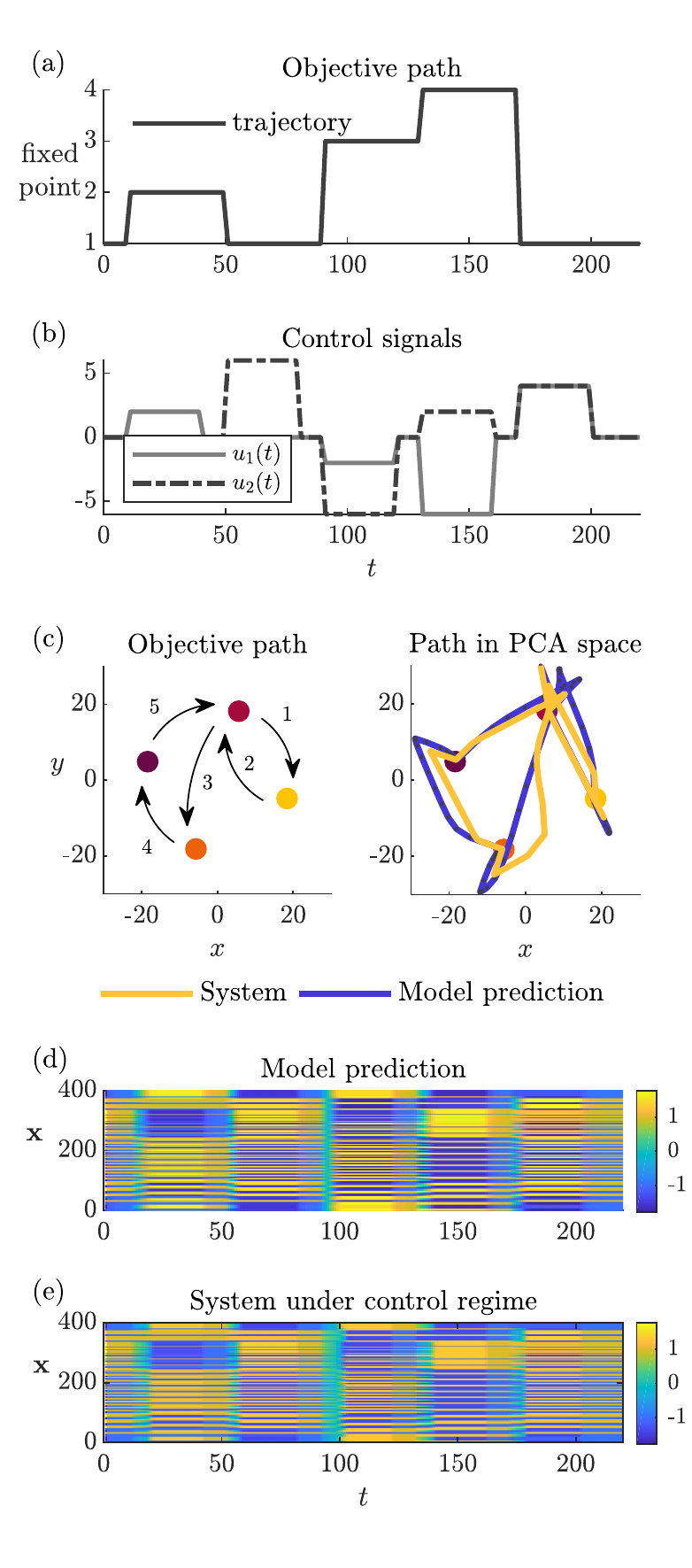}
    \caption{Control example for the hopfield system in Figure~\ref{fig:hop2}.
    (a) Objective path for the system.
    (b) Control signals selected using the systems's stability region maps to move the system between fixed points.
    (c) Objective path in PCA space.
    (d) Predicted and actualized system paths in PCA space.
    (e-f) Predicted and actualized system activity in the original high-dimensional space.}
    \label{fig:hop2_controlled}
\end{figure}

Figure~\ref{fig:hop2} shows timeseries data collected from a Hopfield network storing four memories, the trajectories in PCA space, and the SINDy model approximating the system. The Heaviside function used in the original definition of the model is replaced by a $\tanh$ function in order to make the dynamics smooth and therefore the Taylor series expansion discoverable by the SINDy algorithm. The rapid decay in singular values shows that the system in PCA space is a good representation of the original data, Fig~\ref{fig:hop2}(c). Figure~\ref{fig:hopfield2_1_bif} shows the saddle-node bifurcation curves $C^s(t)$ colored by fixed point locations. These curves indicate the stability regions for the system's four memories in control space. Figure~\ref{fig:hop2_controlled} shows the Hopfield network transitioning between fixed points in the manner specified by the objective path. Transient control signals are selected from the stability regions for each transition and applied to the original system through the first two PCA modes.
We are able to move directly from any stable fixed point to any other stable fixed point in the system according to the stability regions.  Furthermore the system does not exhibit any limit cycles, making it an easy system to control with a feed-forward control signal method.

%

%%%%%%%%%%%%%%%%%%%%%%%%%%%%%%%%%
%%% Subsection - 3D %%%%%%%%%%%%%
%%%%%%%%%%%%%%%%%%%%%%%%%%%%%%%%%

\subsection{Systems with Three-dimensional Intrinsic Dynamics}
Only some systems can be adequately described using the first two PCA modes. Many systems have more complex dynamics that require the use of three or more modes. We demonstrate our control procedure on systems with three-dimensional intrinsic dynamics.
%
%%%%%%%%%%%%%%%%%%%%%%%%%%%%%%%%%%%%
%% Subsection - 3D Hopfield %%%%%%%%
%%%%%%%%%%%%%%%%%%%%%%%%%%%%%%%%%%%%
\subsubsection{Three-dimensional Hopfield}
\begin{figure}[t]
    \centering
    \includegraphics[width=0.9\linewidth]{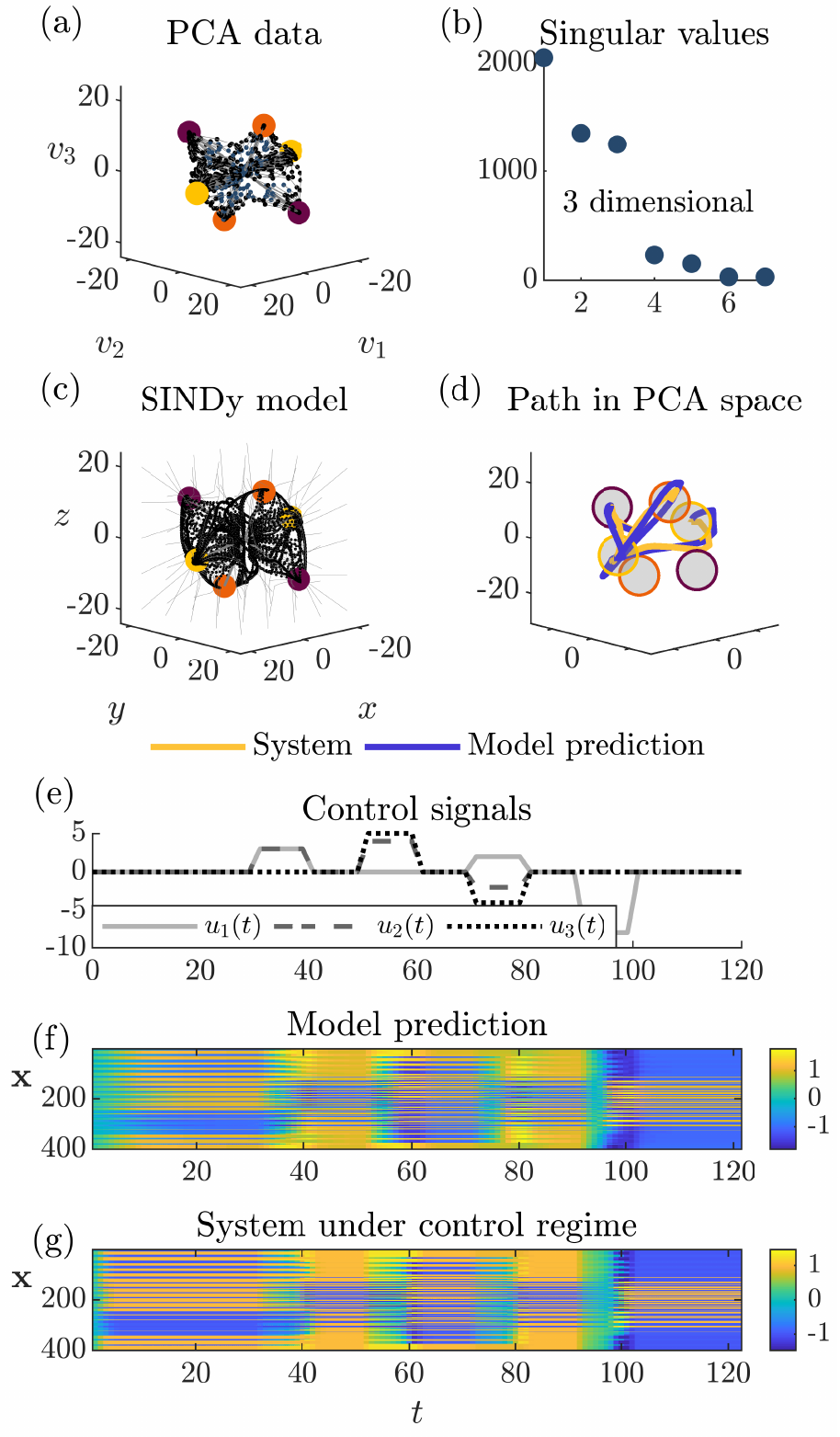}
    \caption{Smooth Hopfield network size $n=400$ with intrinsic 3-dimensional structure. 
    (a) Hopfield data in PCA space. The network converges to one of 6 "memories".
    (b) The first 3 modes capture the majority of the variance in the system.
    (c) SINDy is able to capture a low-dimensional model with similar dynamics that we can use to construct a control regime. 
    (d) Predicted and actualized network activity in PCA space.
    (e) There are now 3 control signals in the control regime as we are using a 3-dimensional system.
    (f-g) Predicted and actualized network activity in the original space under the specified control regime.}
    \label{fig:hop3}
\end{figure}

We inflate the intrinsic dimension of our Hopfield network by increasing the number of memories stored. Figure~\ref{fig:hop3} shows the dynamics of an intrinsically three-dimensional Hopfield network. We use SINDy to fit a three-dimensional nonlinear dynamical system to the data and control the system along three dimensions, $(u_1(t), u_2(t), u_3(t))$. 
While control signals can be selected from analytically derived three-dimensional stability regions, it can be easier in practice to find control signals experimentally by systematically perturbing the SINDy model with control signals to find control regions that will destabilize each fixed point.

%
%%%%%%%%%%%%%%%%%%%%%%%%%%%%%%%%%%%%
%% Subsection - Random 3D %%%%%%%%%%
%%%%%%%%%%%%%%%%%%%%%%%%%%%%%%%%%%%%
\subsubsection{Three-dimensional Random Dynamical System with a Strange Attractor}
\begin{figure}[!ht]
    \centering
    \includegraphics[width=1.05\linewidth]{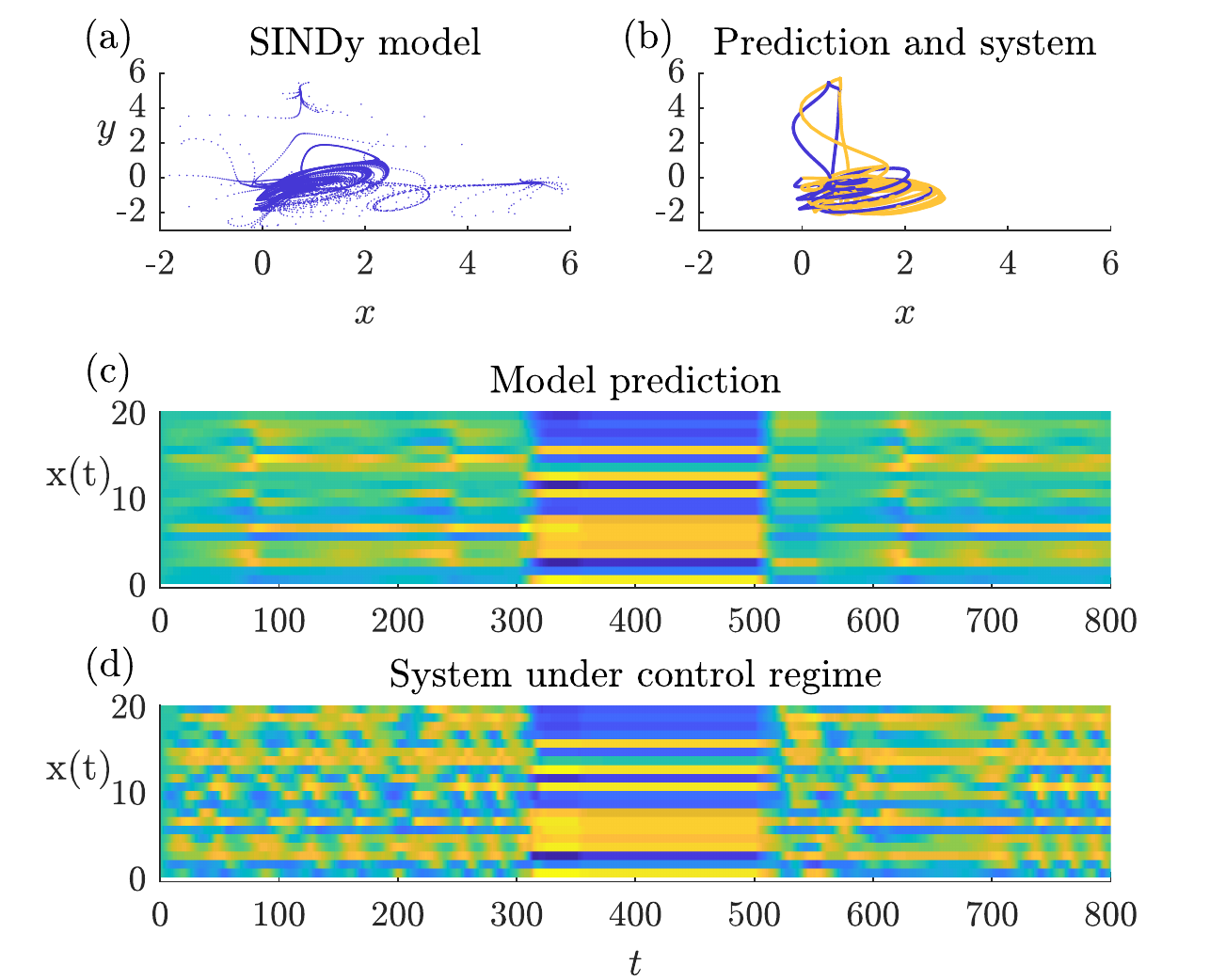}
    \caption{Randomly generated high-dimension system $n = 20$ with a strange attractor. (a) SINDy model with 3 variables captures the stable fixed points and strange attractor in the system.
    (b) Predicted and actualized network activity in PCA space under control.
    (c) Predicted and actualized network activity in the high-dimensional space.
    Notice that although the SINDy model captured the location of the strange attractor, it does not adequately capture its frequency as the model predicts a much slower oscillation that the original system actually produces.}
    \label{fig:random3}
\end{figure}

We generate a random dynamical system that has intrinsic three-dimensional dynamics and use control signals to move the system between a fixed point and a strange attractor in the system. Figure~\ref{fig:random3}(a) shows the SINDy model approximation of the system which captures the general location and stability of the strange attractor and two stable fixed points. While the location of the fixed point and attractor appear to be adequately captured by the SINDy model, Fig.~\ref{fig:random3}(b), the model predicts an oscillation frequency for the strange attractor that is much slower than that observed in the original system, Fig.~\ref{fig:random3}(c-d), meaning that while the model may be able to predict whether the system's state is on or off the strange attractor it cannot predict its location along the attractor. Stability regions for the fixed points can be analytically computed; however, the stability region for the strange attractor cannot be analytically computed and therefore must be experimentally determined by measuring the effects of control perturbations to the SINDy model.

%%%%%%%%%%%%%%%%%%%%%%%%%%%%%%%%%%%
%%%% Discussion %%%%%%%%%%%%%%%%%%%
%%%%%%%%%%%%%%%%%%%%%%%%%%%%%%%%%%%
\section{Discussion}
Our control procedure shifts nonlinear network dynamical systems between attractors using feed-forward control signals derived from the system equations. 
This procedure works for high-dimensional systems that exhibit low dimensional dynamics by developing a control procedure for the system in a reduced space and then projecting the derived control signals back to the original space.
\par
Our framework builds off of previous work in which we used a feed-forward control model to demonstrate how the nematode \textit{C. elegans} may switch between short-term behavioral states using externally generated control signals \cite{morrison_nonlinear_2020}. We observed that the modulation of internal system parameters was another form of control that the nematode could use to implement long-term behavior changes.

\par
While many traditional control strategies used in a wide range of fields can be equivalently viewed as actuating forces applied to complex systems, some control strategies involve modifying internal system parameters or using system feedback, meaning that our control framework does not account for all ways control can be achieved in a network.

\par

In pharmacology \cite{kenakin_pharmacology_2004}, certain agonists can be modeled as feed-forward control signals applied to dynamics on a biochemical network while antagonists may be better represented as shifting network parameters. In clinical neuroscience, deep brain stimulation is used as a form of external control to correct motor output in patients with Parkinson's disease and essential tremor as well as to reduce pain in patients with chronic pain \cite{benabid_deep_2003, perlmutter_deep_2006, kringelbach_deep_2007}. In an alternative setting, deep brain stimulation is used to alter brain structure in patients with treatment-resistant depression and traumatic brain injury \cite{mayberg_deep_2005, schiff_behavioural_2007}. This same medical intervention has the ability to both act as a transient feed-forward control signal to the neural system as well as alter internal structures highlighting how a single control signal applied to a system can have diverse effects.
\par

Network control is also widely implemented in the political and social sciences. Governments desire to optimize economic performance through monetary and fiscal policy \cite{friedman_role_1968, friedman_monetary_1971, arrow_public_1970}. The central bank aims to achieve economic goals by adjusting interest rates, which modify internal economic parameters.  Meanwhile, the federal government pursues these same goals through taxation and government spending which are more easily viewed as feed-forward control strategies.  States often provide external support to insurgent movements in order to further their own interests \cite{byman_trends_2001}; these sponsorships can have a variety of effects on the political network including increasing the power of particular groups or modifying the way groups interact.
During an epidemic, public health experts have various ways of suppressing the spread of an infectious disease. They can achieve public health goals by modifying people's social interactions, which changes interaction parameters in infectious disease models, or by introducing external forces such as vaccinating susceptible individuals or quarantining infected ones \cite{kermack_contribution_1927, hethcote_mathematics_2000}. Given the limited avenues of control and leverage available to health care professionals, policy makers, and organization leaders, all possible forms of control are often implemented concurrently. While this has the potential to maximize influence on the system in question, it also makes identifying the effects of multiple control strategies difficult to differentiate.

Sometimes the systemic changes necessary to bring about a desired outcome are either unknown or not possible to implement in a given network while the use of external control presents a clear path to success. We present a feed-forward control strategy as an alternative form of control for nonlinear network systems that are unamenable to feedback control or alterations to system parameters. 

One limitation of this control strategy is that the locations and stabilities of limit cycles and strange attractors in the system cannot be analytically determined, particularly for systems with more than two dimensions. This limits our ability to assert control guarantees and means that the stability regions of some attractors must be determined experimentally as they cannot be determined analytically. We aim to extend this work by defining the locations and stabilities of limit cycles and attractors under the influence of feed-forward control by using data-driven discovery methods on Poincare maps as presented in \cite{bramburger_poincare_2019}. 

Another limitation of this method is that it can only be applied to high-dimensional systems that can be linearly dimension reduced. Many high-dimensional systems exhibit low-dimensional activity but must be reduced using a nonlinear dimension reduction method. We propose extending feed-forward control to such systems by using autoencoders to perform the nonlinear dimensionality reduction \cite{champion_data-driven_2019} and then transforming the low-dimensional control signals using the autoencoder. Feed-forward control of nonlinear systems lends itself to many developments including optimization, extensions to higher dimensions, and hybridization with feedback control methods.

%
%%%%%%%%%%%%%%%%%%%%%%%%%%%%%%%%%%%
%%%% Conclusion %%%%%%%%%%%%%%%%%%%
%%%%%%%%%%%%%%%%%%%%%%%%%%%%%%%%%%%
\section{Conclusion}
Nonlinear network dynamical systems in the  physical, engineering, biological and social sciences are typically difficult to characterize and control. We develop a feed-forward control procedure for low-dimensional nonlinear systems by deriving control signals as a function of local bifurcations in the system. We extend this method to high-dimensional systems with unknown dynamics by using dimensionality reduction in conjunction with the model-discovery SINDy algorithm. Our technique is demonstrated on canonical and random network dynamical systems of different dimensions with known and unknown dynamics. We propose this method of control as an alternative to feedback control and as a framework for understanding how actuating forces can be used in nonlinear systems to shift systems between stable states.

To our knowledge, our method provides a principled mathematical architecture that integrates dimensionality reduction, bifurcation theory, and data-driven discovery of dynamics to construct a feed-forward control model for nonlinear, networked dynamical systems.  Our feed-forward control techniques can be interpreted and used to regulate the dynamics of networked dynamical systems by exploiting the dominant, low-dimensional subspaces on which the dynamics evolves. Our mathematical architecture generates a set of actuation signals that, when applied, are able to control the original high-dimensional system. Using bifurcation theory, we find  collections of feed-forward control signals that will force convergence to desired objective states, allowing us to move the system from one fixed point of the system to another in a principled manner. Specifially, we can destabilize a given fixed point by making it undergo a saddle node or Hopf bifurcation, while simultaneously making the target fixed point an attractor.  This creates a pathway with the feed-forward signals from one fixed point to another.  

The potential applications of this control framework  are numerous.  From neuroscience to powergrids, networked dynamical systems appear in almost every branch of quantitative study.  Many modern systems of study are high-dimensional and characterized by nonlinear dynamics, thus they require new mathematical methods to characterize their behavior and exploit their intrinsic dynamics.  The methods developed here rely on the simple observation that most high-dimensional dynamical systems manifest low-dimensional dynamics.  Thus the low-dimensional subspaces on which the system evolves can be exploited for control, as is done in this manuscript.  The number of application areas for which this is true is quite diverse and extensive, including neuroscience~\cite{kunert_low-dimensional_2014,kato_global_2015,fieseler_unsupervised_2020,jones_natural_2007,rabinovich_transient_2008,shlizerman_data-driven_2014,delahunt_biological_2018}, powergrids~\cite{dylewsky_engineering_2019}, disease modeling~\cite{kermack_contribution_1927, hethcote_mathematics_2000} and social networks~\cite{morrison_community_2019}.

% use section* for acknowledgment
\section*{Acknowledgment}
This work was funded in part by the Big Data for Genomics and Neuroscience Training Grant under the National Institute of Health Grant No. 5T32LM012419-04 and by the Army Research Office under Grant No. W911NF1910291.
JNK also acknowledges support from Air Force Office of Scientific Research (AFOSR) grant FA9550-
17-1-0329.

% Can use something like this to put references on a page
% by themselves when using endfloat and the captionsoff option.
\ifCLASSOPTIONcaptionsoff
  \newpage
\fi

\bibliographystyle{IEEEtran}  
\bibliography{me,Nathan,Nonlinear_control}%

% that's all folks
\end{document}